\documentclass[trsc]{paper} 

\DoubleSpacedXI

\usepackage[utf8]{inputenc}
\usepackage{pgfplots}
\DeclareUnicodeCharacter{2212}{−}
\usepgfplotslibrary{groupplots,dateplot}
\usetikzlibrary{patterns,shapes.arrows}
\pgfplotsset{compat=newest}

\usepackage{bm}
\usepackage{caption}
\usepackage{subcaption}
\usepackage{rotating}
\usepackage{mathtools}
\usepackage{pgfplots}
\usepackage{algorithm}
\usepackage{algpseudocode}
\usepackage{booktabs}
\usepackage{cases}
\usepackage{enumerate}
\usepackage{hyperref}

\pgfplotsset{compat=1.17}

\usepackage{xcolor}

\usepackage{tikz}
\usetikzlibrary{arrows}
\usetikzlibrary{shapes}
\definecolor{tempcolor}{RGB}{244,244,244}
\usetikzlibrary{positioning}
\tikzset{
    split/.style = {shape=rectangle,
                     draw, align=center,
                     fill=tempcolor},
    clust/.style = {align=center,
                    draw,
                    rectangle}
                     }

\tikzstyle{process} = [rectangle, minimum width=3cm, minimum height=1cm, text centered, draw=black, fill=orange!30]
\tikzstyle{arrow} = [thick,->,>=stealth]

\usepackage{pgfplots}
\usetikzlibrary{
    patterns,
}
    \pgfplotsset{
        compat=1.7,
        my ybar legend/.style={
            legend image code/.code={
                \draw [##1] (0cm,-0.6ex) rectangle +(2em,1.5ex);
            },
        },
    }

\algblockdefx{FORP}{ENDFORP}[1]%
  {\textbf{for }#1 \textbf{do in parallel}}%
  {\textbf{end for}}
  
\usepackage{todonotes}

\RUNAUTHOR{ \!} 
\RUNTITLE{A Data-Driven Methodology for Last-Mile Routing} 
\definecolor{elegantblue}{RGB}{18, 52, 89}

\newcommand{\tc}[1]{{\textcolor{black}{#1}}}

\usepackage{hyperref}
\hypersetup{
    colorlinks=true,
    linkcolor=blue,
    filecolor=blue,  
    citecolor=elegantblue,
    urlcolor=elegantblue,
}
\usepackage{natbib}
 \bibpunct[, ]{(}{)}{,}{a}{}{,}%
 \def\bibfont{\small}%
\TheoremsNumberedThrough     

\EquationsNumberedThrough    

\begin{document}




\TITLE{Learn Global and Optimize Local: A Data-Driven Methodology for Last-Mile Routing} 

\ARTICLEAUTHORS{%
\AUTHOR{Mayukh Ghosh$^{1}$, Alex Kuiper$^{1}$, Roshan Mahes$^{1,2}$, Donato Maragno$^{1}$}
\AFF{1 Amsterdam Business School, University of Amsterdam, 1018 TV Amsterdam, the Netherlands \\
2 Korteweg-de Vries Institute for Mathematics, University of Amsterdam, 1098 XG Amsterdam, the Netherlands \\ 
\EMAIL{m.ghosh@uva.nl}; \EMAIL{a.v.mahes@uva.nl}; \EMAIL{d.maragno@uva.nl}; \EMAIL{a.kuiper@uva.nl}}
}

\ABSTRACT{

In last-mile routing, the task of finding a route is often framed as a Traveling Salesman Problem to minimize travel time and associated cost. However, solutions stemming from this approach do not match the realized paths as drivers deviate due to navigational considerations and preferences. To prescribe routes that incorporate this tacit knowledge, a data-driven model is proposed that aligns well with the hierarchical structure of delivery data wherein each stop belongs to a zone---a geographical area. First, on the global level, a zone sequence is established as a result of a minimization over a cost matrix which is a weighted combination of historical information and distances (travel times) between zones. Subsequently, within zones, sequences of stops are determined, such that, integrated with the predetermined zone sequence, a full solution is obtained.

The methodology is particularly promising as it propels itself within the top-tier of submissions to the \textit{Last-Mile Routing Research Challenge}, while it maintains an elegant decomposition that ensures a feasible implementation into practice. The concurrence between prescribed and realized routes 
underpins \tc{the adequateness of a hierarchical breakdown of the problem, and the fact that drivers make a series of locally optimal decisions when navigating. Furthermore, experimenting with the balance between historical information and distance exposes that historic information is pivotal in deciding a starting zone of a route. The experiments also reveal that at the end of a route, historical information can best be discarded, making the time it takes to return to the station the primary concern.}




}%


\KEYWORDS{last-mile logistics, route prediction, navigation, traveling salesman problem, data-driven optimization}

\maketitle

%

\newpage

\section{Introduction}

Last-mile routing is often the shortest, yet most complex and expensive leg in distribution systems. It accounts for about 28\% of the total costs of transportation from the distribution center to the final customer \citep{goodman2005whatever,gevaers2011characteristics}. Last-mile operations also have a large impact on the livability, because of traffic congestion and pollution it creates~\citep{taniguchi2002modeling,chopra2003designing}. This pressure is set to increase, because the trend of urbanization continues; around 68\% of the global population will reside in urban areas by 2050, according to a report by the \cite{UN2018}. \tc{Currently e-commerce sales are  estimated at about \$709.78 billion, which is projected to increase as the coronavirus has  accelerated a channel shift to e-commerce, see \cite{eMarketer2020}. The same trends are echoed by a recent report of the World Economic Forum \citep{deloison2020future}, which further points out that also the number of same-day deliveries will grow due to customers ever-increasing expectations, amplifying the importance of last-mile operations and its impact on the ecosystem.} 

\tc{Considering parcel delivery in last-mile operations,} the main goal is to ensure that products arrive in good shape. In addition, it should be on time~\citep{gunasekaran2004framework}, i.e., the time windows that have been communicated will be met. As routes require a multitude of stops to be visited, it is critical that a proposed route needs to be followed up as deviations can incur redeliveries or dissatisfied customers (by not meeting time slots), unforeseen time losses along the route and additional operational expenses~\citep{boyer2009last,GEVAERS2014398}.

Despite the benefits of adhering to a proposed `optimal' route are clear, in practice---reported by managers in last-mile delivery (e.g., Amazon)---drivers often deviate from proposed  routes. This discrepancy between theoretical route planning and real-life route execution is an important gap in the literature on last-mile routing, and has remained an open research question, because of the typical inaccessibility of real-life data. Recently, Amazon in combination with MIT's Center for Transportation \& Logistics hosted a challenge \citep{routingchallenge}, wherein by means of offering a vast routing dataset, researchers around the world have been invited to bridge this gap. In this work we answer that call by providing a complete methodology that proposes routes that adequately reflect drivers' tacit knowledge and preferences by learning from historic route data. 


In the US, Amazon steadily holds the number one position as the largest e-commerce retailer having a market share of 38.0\%~\citep{eMarketer2020}, \tc{and inherently is responsible for the majority of traffic associated to last-mile delivery. Therefore, by applying the methodology on the data associated to Amazon will set an example for many other organizations facing the bewildering setting of last-mile delivery. Although we consider it as a predictive problem, which helps uncovering on a conceptual level the driver's decision process when determining a route, we emphasize that the approach can also be used to prescribe routes, which are more feasible from a driver's perspective.}


\begin{figure}
    \centering
    \includegraphics[width=\textwidth]{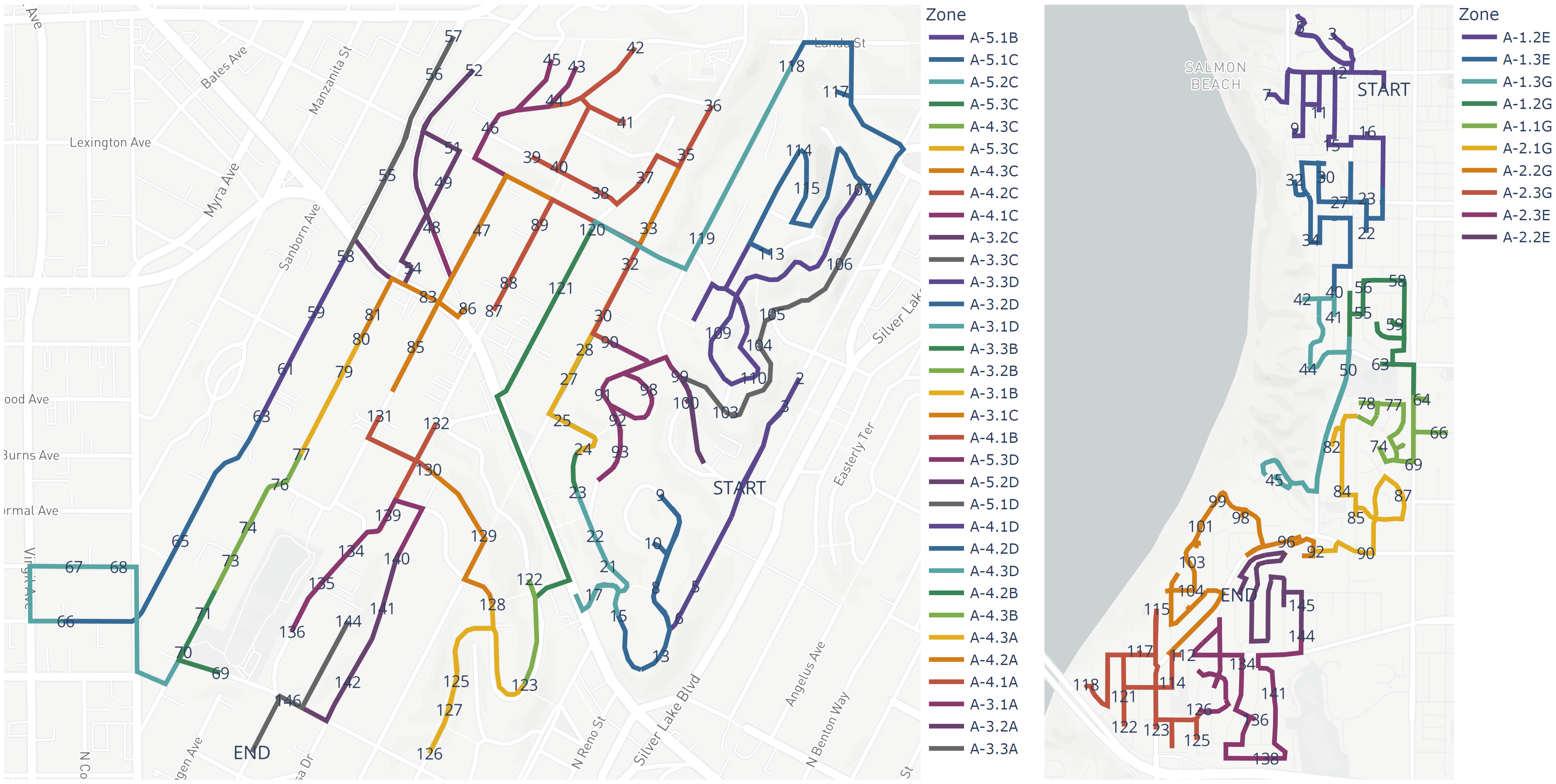}
    \caption{Examples of two realized routes by using \cite{OpenStreetMap}. The numbers represent the order in which the stops are visited and the colors represent the different zones within the routes.}
    \label{fig:example_route}
\end{figure}

To consider the problem at a more detailed level, two delivery drivers' routes are given in Figure~\ref{fig:example_route} as an increasing sequence of numbers; each number represents a stop. Note that not all stops are shown as for readability purposes.  Along the route the driver traverses through different zones, given by their zone ids, which are (preset) \textit{geographical areas}. Comparing the routes in the left (more rectangular grid) and the right map, it is observed that the size of a zone ranges from a part of a street to a set of multiple streets. 

\begin{figure}
\centering
\begin{subfigure}{.49\textwidth}
\centering
\includegraphics[width=\textwidth]{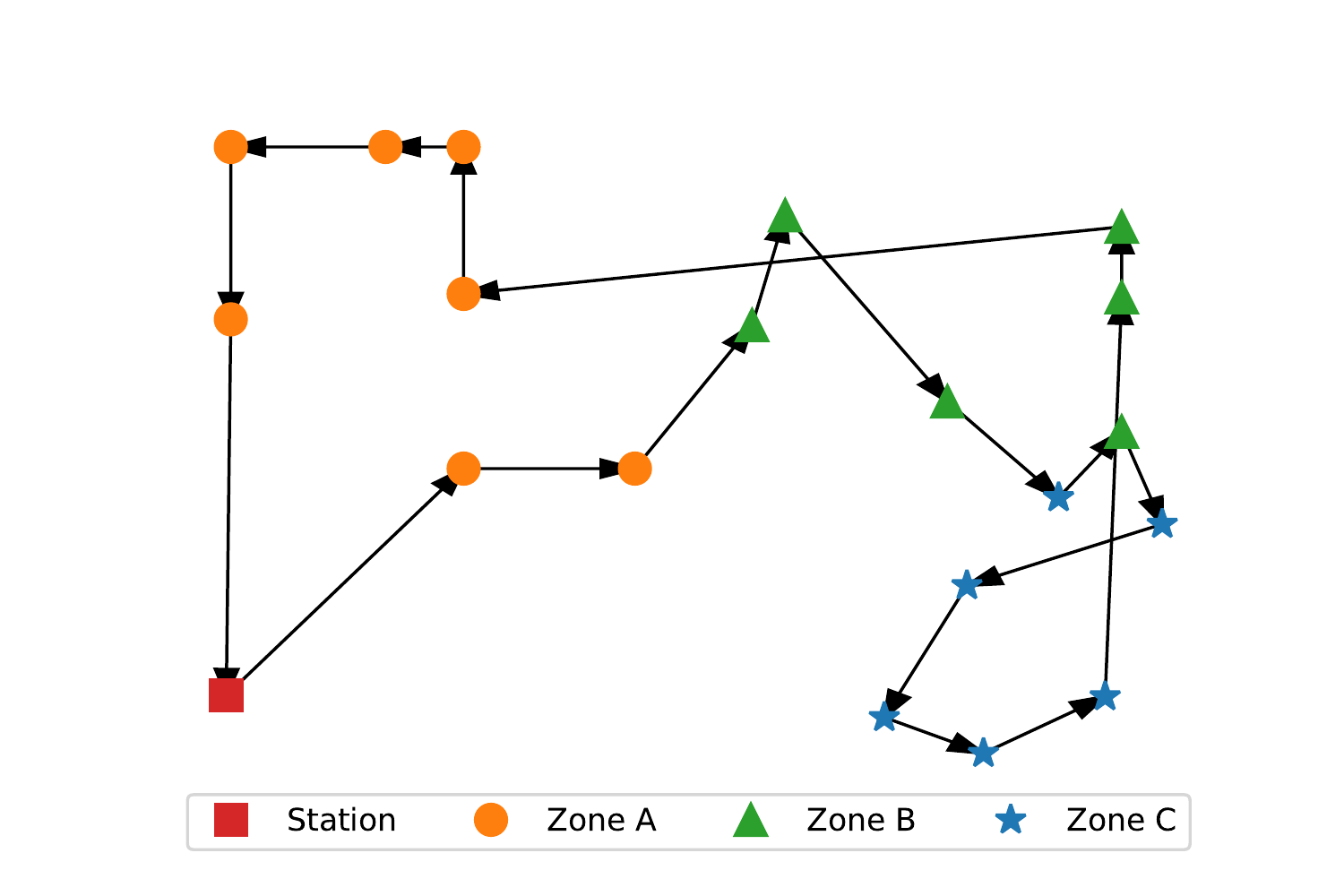}
\caption{Nearest Neighbor solution.}
\label{fig:nn_solution}
\end{subfigure}
\begin{subfigure}{.49\textwidth}
\centering
\includegraphics[width=\textwidth]{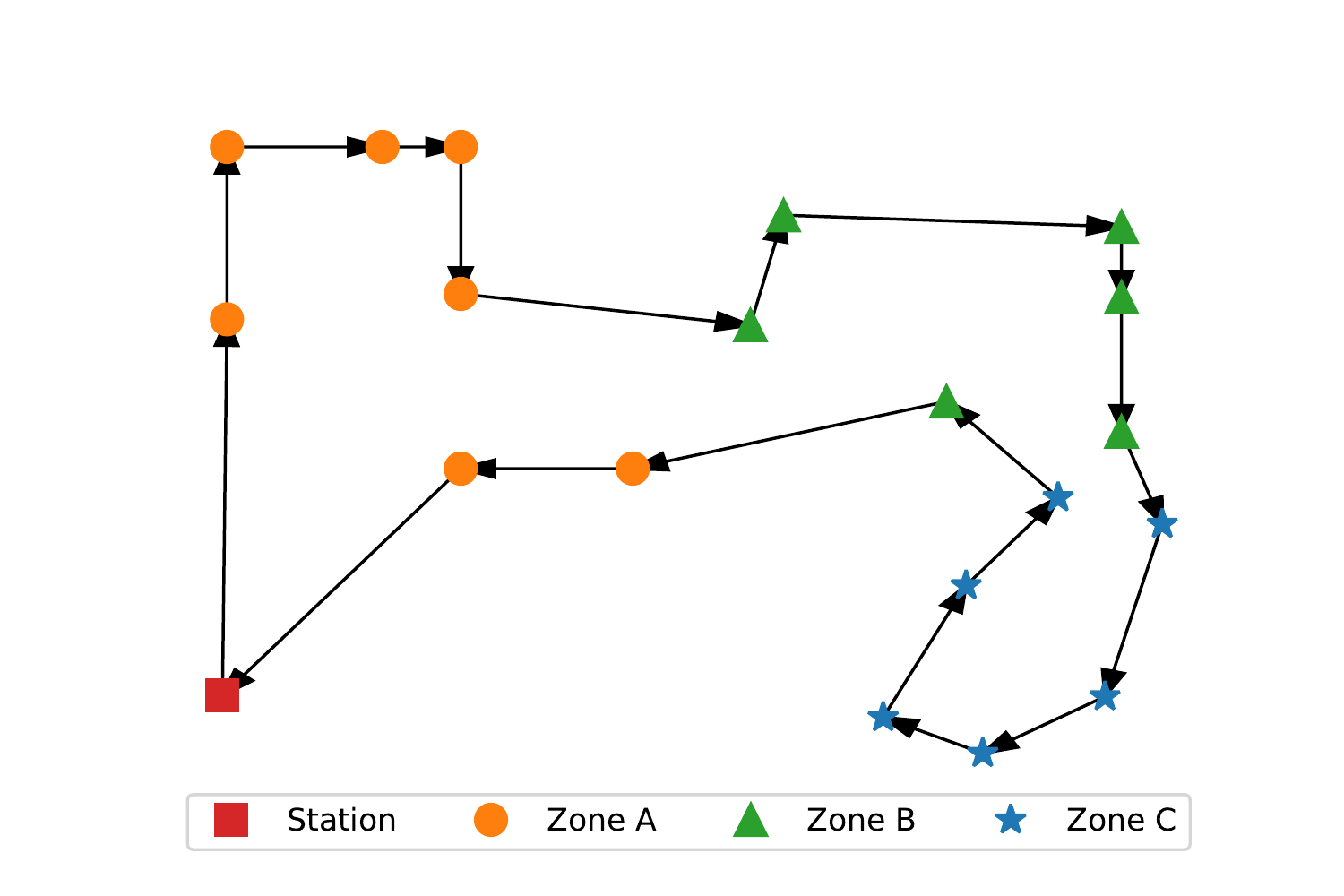}
\caption{TSP solution.}
\label{fig:tsp_solution}
\end{subfigure}
\caption{An example of a heuristic (a) and the solution when using an optimization model (b).}
\label{fig:nn_vs_tsp}
\end{figure}

Leveraging historical routes to learn recurrent patterns is feasible on a zone level rather than at a stop level. This is due to the fact that a zone appears in multiple routes while a stop does not (or rarely). The intrinsic importance of a zone sequence depends on how the zone is established. Ideally, a zone should be constructed such that every stop in it is visited before moving to the next zone (as for example in Figure~\ref{fig:example_route}). However, with this presumption, classical heuristics and optimization models fall short to exploit this information. Two classical examples are illustrated in Figure~\ref{fig:nn_vs_tsp}. On the left, we find the route obtained by applying the nearest neighbor heuristic, which is sometimes argued to mimic the driver’s behavior, while on the right the route that minimizes distance, as the solution of the corresponding Traveling Salesman Problem (TSP).

On top of the condition that drivers visit all stops in a zone before moving to another, the challenge is to determine the sequence of zones a driver will follow. In Figure~\ref{fig:nn_vs_tsp}, there are six possible permutations of zone sequences, but in real-life settings (Figure~\ref{fig:example_route}) there are easily 30 zones leading to 30! unique zone sequences. The sequence might be based on information at ground level, navigational considerations, or merely a preference, but most importantly an approach should be able to learn this tacit information from realized routes.

 
 
The contributed methodology divides the problem into two parts: on the zone level, it aims to infer the sequence of zones after which, within each zone, the sequence of stops is determined. \tc{It stems from the idea that first a routing problem is abstractified into zones,} for which a sequence is established by solving a TSP formulation that balances historical zone transition information, resembling the tacit knowledge, and distances. So the methodology learns on a global level, which aligns with the intuition that humans congregate route sequences on a higher level than remembering exact stop sequences. Secondly, the problem in each zone is of relatively low complexity and thus a driver is capable to determine an optimal sequence out of the set of stops that minimizes distance or travel time. In the end, the performance of the methodology confirms the competency of the approach. 

\tc{In the next section, we provide the literature review, which outlines important elements to our methodology. Section~\ref{sec:model} describes the methodology, and is split into learning and prediction. In Section~\ref{sec:perf}, we introduce the dataset on which the methodology will be applied and provide the route prediction metrics aligned with the challenge. Using the case study data we carry out various experiments and demonstrate the methodology. Furthermore, we compare the performance to the aforementioned benchmarks as well as to the submissions of the competition. In Section~\ref{sec:concl}, we present our concluding remarks and implications to practice, followed by a discussion and some suggestions for future research.}

\section{Literature Review}\label{sec:literaturereview}
Solving a last-mile routing problem can be translated to a TSP, i.e., a set of stops should be visited---starting from and ending at the station. However, in this research the goal is to predict a driver's route using historic data. Therefore, in this pursuit, elements of different disciplines are combined. More specifically, we  consider classical TSP optimization frameworks, hierarchical structures as an approximation (such as cluster first, routing second), human approaches to solving TSPs, and lastly, learning from route data. In relation to our work, we elaborate further on these themes below. 

\subsection{Traveling Salesman Problem}\label{TSP}
The TSP with Euclidean distance is notorious for being hard and is proven to be NP-complete \citep{papadimitriou1977euclidean}, but for finding the optimal or near-optimal tour there is a rich literature~\citep{lawler1985traveling,applegate2011traveling}. Also, several variations of the standard TSP have been introduced in literature, to name a few: an open TSP (OTSP) where the driver does not necessarily return to the depot; TSP with time windows to meet~\citep{gendreau1998generalized}; additionally allowing pick-ups such that the capacity of the truck is not exceeded~\citep{gendreau1999heuristics}; or the covering Salesman problem (CSP)~\citep{CURSHI89}, where a minimal tour is found such that all stops are in the vicinity of a stop on the tour (e.g., dropping goods at neighbors). 

Acquainted to the TSP is its generalization known as \textit{Vehicle Routing Problem (VRP)}; see \cite{Dantzig1959} for its formulation and \cite{braysy2005vehicle,MONTOYATORRES2015115} for recent surveys on multiple depots VRPs with variants that include time windows, split delivery, heterogeneous fleet, periodic deliveries, and pick-up. Also, various studies have considered the ecological footprint as presented in the literature review given by  \cite{Canhong2014}. Finally, also for the more general VRP, numerous heuristics and metaheuristics have been proposed over the last decades \citep{braysy2005vehicle2,Speranza2014, Roneeta2020}. 

Both TSP and VRP  have seen many different integer programming formulations \citep{Orman2007ASO}. In our model, we adopt the Miller-Tucker-Zemlin (MTZ) formulation to solve the TSP, introduced for the first time by \cite{Miller1960IntegerPF}. The MTZ formulation has the advantage of being intuitive and straightforward to implement and works well when the number of stops (variables) to consider is relatively small, fitting our instances. As a downside, it has been proven to lead to a weak relaxation, leading to higher computation times, see e.g., \cite{campuzano2020accelerating} for a discussion. In our case, as we break down the full problem into smaller ones, this is not an issue, but as desired any other TSP formulation can be adopted. 
 
\subsection{Hierarchical Structure}
To reduce the computational complexity, the problem can be broken down into smaller instances, which are computationally much less involved, and connect these together. For example, in the seminal work by \cite{karp1977probabilistic}, the proposed partitioning algorithms are asymptotically optimal heuristics in the case of uniformly distributed stop locations, i.e., the error tends to zero with probability one as the length of a route (i.e., the number of stops) increases.

Both \cite{8404041} and \cite{jiang2014hierarchical} decompose the problem with small-scale nodes by relying on clustering algorithms. Each subproblem is subsequently optimized and the center nodes of the subproblems constitute a TSP in itself. Connecting all local tours in the order of the upper layer problem generates approximative solutions with significantly reduced computation times. Such an approach is also found from an empirical point of view, see \cite{vickers2003roles}.
 
 Moreover, evidence to support the use of a hierarchical structure in practice is found in~\cite{graham2000traveling}. They show that solution methods originating from artificial intelligence or operations research algorithms are insufficiently capable to mimic the human approach of solving TSPs. They introduce a hierarchical model by means of a pyramid algorithm on the visual representation. This algorithm is capable to render human-alike solutions to the various TSP instances where classical algorithms fail in this task, see also~\cite{pizlo2006traveling} for a more refined algorithm. More importantly, the aforementioned works hint at the use of a hierarchical breakdown, because of lower computational complexity and being concurrent with practice.

\subsection{Human Navigation}\label{sec:human_nav}
The TSP lends itself to a broad range of experimental studies as its goal of minimization of the route is easy to understand and visualize. Many experimental studies have been devoted to visual versions of the TSP~\citep{macgregor1996human,van2003convex,vickers2003roles}. \cite{macgregor1996human} show that humans outperform well-known TSP heuristics and are in small problem sizes capable to be in 1\% from optimality. 

Humans rely on various tactics to generate near-optimal solutions for the TSP. The potential of these tactics is further demonstrated by the fact that the time needed increases only in a linear fashion compared to the problem size~\citep{pizlo2006traveling,dry2006human}. One of these tactics is to consider the problem first globally, considering the tour that visits all `exterior' points, to which interior points are inserted; this resembles the so-called convex hull approach~\citep{macgregor1996human,vickers2003roles}. 
Argued as the underlying motivation for the convex hull approach, is the avoidance of crossings in a tour \citep{van2003convex}. The motivation behind this tactic is the intuition that a cross in a tour is suboptimal, which is even a fact for metric TSPs. In an OTSP, where the solution is not a tour, the starting point can have a profound impact on the performance, see \cite{sengupta2018planning}. However, they also report that humans are quite capable to select a `right' starting point.

\cite{wiener2009planning} find in a series of experiments that, when navigating, a coarse route is stipulated first. This route visits a set of `regions' after which it is `optimized' on a detailed level along the way. Furthermore, if the problem size increases, the problem is divided into more regions to make the problem manageable and approachable. Such a global-to-local approach echoes the hierarchical decomposition of first determining the zone sequence on a global level, and subsequently a series of local problems to determine the sequence in which stops are to be visited in each zone. Additionally, another experiment from~\cite{wiener2003fine} reveals that a segmentation into zones affects the route planning and navigational behavior as it primes a driver to approach the problem from such a structure. \tc{A similar effect can be expected when drivers are presented a route, for example a solution from routing software as in the Amazon challenge, from which we take the data to demonstrate our methodology~\citep{routingchallenge}.} 

\subsection{Route Prediction}\label{sec:learningTSP}
Typically, the problem of finding a sequence of stops is positioned in the field of combinatorial optimization, but in this research the goal of find an optimal sequence to minimize an objective function does not take center stage. Here, the goal is how route data can be used to predict routes that drivers (will) follow; for an example where learning takes place in order to find the optimal route, see~\cite{kool2019attention}. 

The capability to learn and predict (part of) the routes is demonstrated in~\cite{krumm2008markov}. The model is trained from drivers' long-term trip history using GPS data. It uses a Markovian approach, i.e., on the basis of the last road segment, which can also be adjusted to incorporate more history, it predicts the next segment the driver will take. 
In the same vein, \cite{ye2015method} propose a route prediction method based on a hidden Markov model that can accurately predict an entire route early in the trip. \cite{wang2015building} also employ a Markov model; their algorithm relies on a probability transition matrix that is developed to represent the knowledge of the driver’s preferred links and routes. For the VRP, \cite{canoy2019vehicle} show the potential of using a Markov model in an optimization framework. By constructing a transition probability matrix, based on historical data, and by exploiting the VRP structure, they render solutions that resemble actual route plans much better than relying on a distance metric.
 
The previous works demonstrate that a Markov model is a powerful tool to learn from historical data and to use these in predictions. We adopt such a model by using historical information, which is weighted against a distance metric. The full methodology is detailed in the next section. 

\section{Methodology}\label{sec:model}

The methodology that we propose relies on the fourfold of elements reviewed in Section~\ref{sec:literaturereview}. We impose a hierarchical structure on two levels: finding the zone sequence on a global level and finding the stop sequence within each zone locally. Such an approach is in line with evidence from human navigation, and more importantly matches the structure of typical routing data---each stop belongs to a zone. We learn on the global level the preference to traverse from one zone to another by adopting a Markov model, which we combine with the distances between zones to form the cost matrix for the global problem. Relying on standard methods, we solve the TSP on the global level to determine the zone sequence, and a series of TSPs on the local level, which mimics the driver capabilities to find near-optimal solutions in small instances. Finally, by patching the local solutions, we deliver a full prediction of the route. 

To learn the drivers' preference over zones we use the fact that each route starts at a station $s$---also commonly known as depots, where $s$ belongs to the set of stations $\mathcal{S}$. Then for each stop in a route, the data consists of the \textit{GPS coordinates}, which can be used to compute distances between stops, but the methodology can also be used when the costs of traversing from one stop to another are given, e.g., (expected) travel times. Furthermore, each stop within the route should be given an identifier, which relates the stop to belong to a geographical area, i.e., a zone. The zones should not overlap and therefore they form a segmentation. As a side note, the zone ids should be consistent within the set of routes in the dataset. Specifically in relation to zones there are two phases, which form the basis of our approach:
\begin{itemize}
    \item \textit{Learning.} Extracting zone sequences from historical data enables the counting and collection of zone transitions in a count matrix. 
    \item \textit{Prediction.} For new routes, a combination of a distance and count matrix is used in a TSP formulation which solution provides the zone sequence; after which within zones the stop sequences are iteratively determined by means of a series of OTSPs.
\end{itemize}
 

\subsection{Learning Zone Preferences}\label{build}
For learning, we exploit the information contained in route data to extract information about drivers' preferences at the level of zones.
Each route consists of an ordered series of $n+1$ stops, namely $n$ delivery nodes (this number differs per route) and the station. Therefore, a route is represented as a sequence of tuples, $\texttt{route}=((x_0,y_0,z_0),\ldots,(x_{n},y_{n},z_{n}))$, with $(x_i,y_i)\in\mathbb{R}^2$ indicating the stop's geographical location, and $z_i$ the corresponding zone id, with $z_i\in\left\{Z_1,\ldots,Z_M\right\}$ for $i\in\{0,\ldots,n\}$. As the data provides the stop sequence per route, we have to convert it to the higher-level sequence of zones.

Given the observations in the introduction and a navigational perspective that all stops within a zone are completed before a driver moves to another zone (Section~\ref{sec:human_nav}), we determine for each route the sequence of unique zones. To illustrate the procedure let us consider a route with $n$ delivery nodes. First, we consider the sequence of zones for each stop $(z_1, \dots, z_n)$. Then, the route’s zone sequence is reduced to the sequence of pairs $((\zeta_1, \tau_1)\ldots, (\zeta_\ell, \tau_\ell))$, where each pair ($\zeta_i, \tau_i$) represents how many nodes ($\tau_i$) belonging to the same zone ($\zeta_i$) are visited in a row. When the $\zeta_i$ are pairwise distinct, we consider $(\zeta_1,\dots,\zeta_\ell)$ to be the desired sequence of unique zones.
\begin{example}
The following route consisting of six stops over three zones is reduced into a sequence of three unique zones:
$$
(z_1, z_2, z_3, z_4, z_5, z_6) = (Z_1, Z_1, Z_1, Z_2, Z_3, Z_3) \rightarrow ((\zeta_1,3),(\zeta_2,1),(\zeta_3,2)) \rightarrow (\zeta_1, \zeta_2, \zeta_3) = (Z_1,Z_2,Z_3). 
$$
\end{example}
In case of a zone's reoccurrence ($\zeta_i=\zeta_j$ for an $i<j$), we only keep the occurrence with the highest number of stops by comparing $\tau_i$ to $\tau_j$, or keep the earliest appearance in case of a tie, here that is $\zeta_i$. This procedure is repeated until we are left with a sequence of singly occurring zones.  Two extreme, yet illustrative, examples of this procedure are given below. 
\begin{example}
The reduction of two routes consisting of six stops over three zones under the procedure outlined above is illustrated below. First, zone $Z_3$ is put last as $\tau_1=1$ and $\tau_4=2$:
\begin{align*}
(Z_3, Z_1, Z_1, Z_2, Z_3, Z_3) &\rightarrow
((\zeta_1,1),(\zeta_2,2),(\zeta_3,1),(\zeta_4,2)) \rightarrow (\zeta_2, \zeta_3, \zeta_4) = (Z_1,Z_2,Z_3).
\end{align*}
In the next case, zone $Z_1$ occurs thrice, but breaking the tie twice puts zone $Z_1$ at its first occurrence in the sequence:
\begin{align*}
(Z_1, Z_3, Z_1, Z_2, Z_2, Z_1) &\rightarrow
    ((\zeta_1,1),(\zeta_2,1),(\zeta_3,1),(\zeta_4,2),(\zeta_5,1))
    \rightarrow (\zeta_1, \zeta_2, \zeta_4) = (Z_1,Z_3,Z_2). 
\end{align*}
\end{example}
Evidently, not all zones are visited within each route, therefore we keep track of all $M$ zones that have been visited at least once in any route. Having reduced each route to its corresponding zone sequence, we compute an asymmetric count matrix $N$ where each entry $N_{ij}$ represents the number of times a driver went from the $i$-th to the $j$-th zone, so that $i$ and $j$ range from $1$ to $\text{$M+S$}$, where $S$ is the number of stations in the dataset. The pseudocode related to the procedures outlined is provided in the Appendix.

\subsection{Predicting the Zone Sequence}\label{sec:pred_zsq}
In the next phase, we predict the route that a driver will take. For that purpose, we 
translate the count matrix, created in the previous phase, to a transition matrix $P$ via $P_{ij}=\dfrac{N_{ij}}{\sum_{j=1}^m N_{ij}}$, so that $P_{ij}$ reflects the Markovian probability of transitioning from zone $i$ to zone $j$. 

Besides the zone sequences, the locations of the station and zones are necessary ingredients in the suspected trade-off a driver makes. For this purpose, we define the center of a zone, if not given, by computing the mean latitude and longitude of all stops that were visited in that zone.  Having the zone centers' location, we can compute the distances $\delta_{ij}$ between zone $i$ and zone $j$. We normalize the distances by dividing each element through the largest, i.e.,  $\max_{i,j}\delta_{ij}$, and store these distances in a matrix \tc{$D \in \mathbb{R}^{(M{+S})\times(M{+S})}$}.

As an alternative to the Euclidean distances in~$D$, a sensible choice is to use the travel times (or another measure for costs between stops) as a link to measure `distance' between zones. However, in the case when we are only given route data, we can proceed with a two-step procedure to have an approximation of the travel times between zones. Since stops are typically assigned to zones, we first identify for each zone the closest stop within the route to the zone center by using the Euclidean distance metric. Second, once each zone has a stop assigned, we extract the travel times between these stops to compose a matrix $T$ which consists of elements $T_{ij}$ corresponding to normalized travel times from station/zone $i$ to $j$, i.e., $T_{ij}=\frac{t_{ij}}{\max_{i,j} t_{ij}}$ where $t_{ij}$ are the provided travel times. We will compare the performances of both options. 

Since we are to minimize a cost matrix, we reverse the transition matrix to $1-P_{ij}$ in the cost matrix to arrive at the following linear combination with $\omega\in[0,1]$ being a weight parameter:
\begin{equation}\label{eq2}
C_{ij} = \omega D_{ij} + (1 - \omega)(1 - P_{ij}),
\end{equation}
where $D_{ij}$ can be interchanged with $T_{ij}$, and additionally we set the diagonal entries $C_{ii} = 0$. In fact this problem can be considered as a TSP with unknown costs as we do not a priori know how distance is evaluated against likelihood of zone-to-zone transitions, which contain the implicit tacit knowledge and drivers' preferences. So a good value of $\omega$ is still to be determined. Next, we are able to use this cost matrix in a TSP formulation, which optimal solution constitutes a tour through all zones in which one or more stops are located. For details about the implementation we refer the reader to the Appendix.



\subsection{Predicting the Stop Sequence}
\begin{figure}
\begin{subfigure}{.49\textwidth}
\centering
\includegraphics[width=\textwidth]{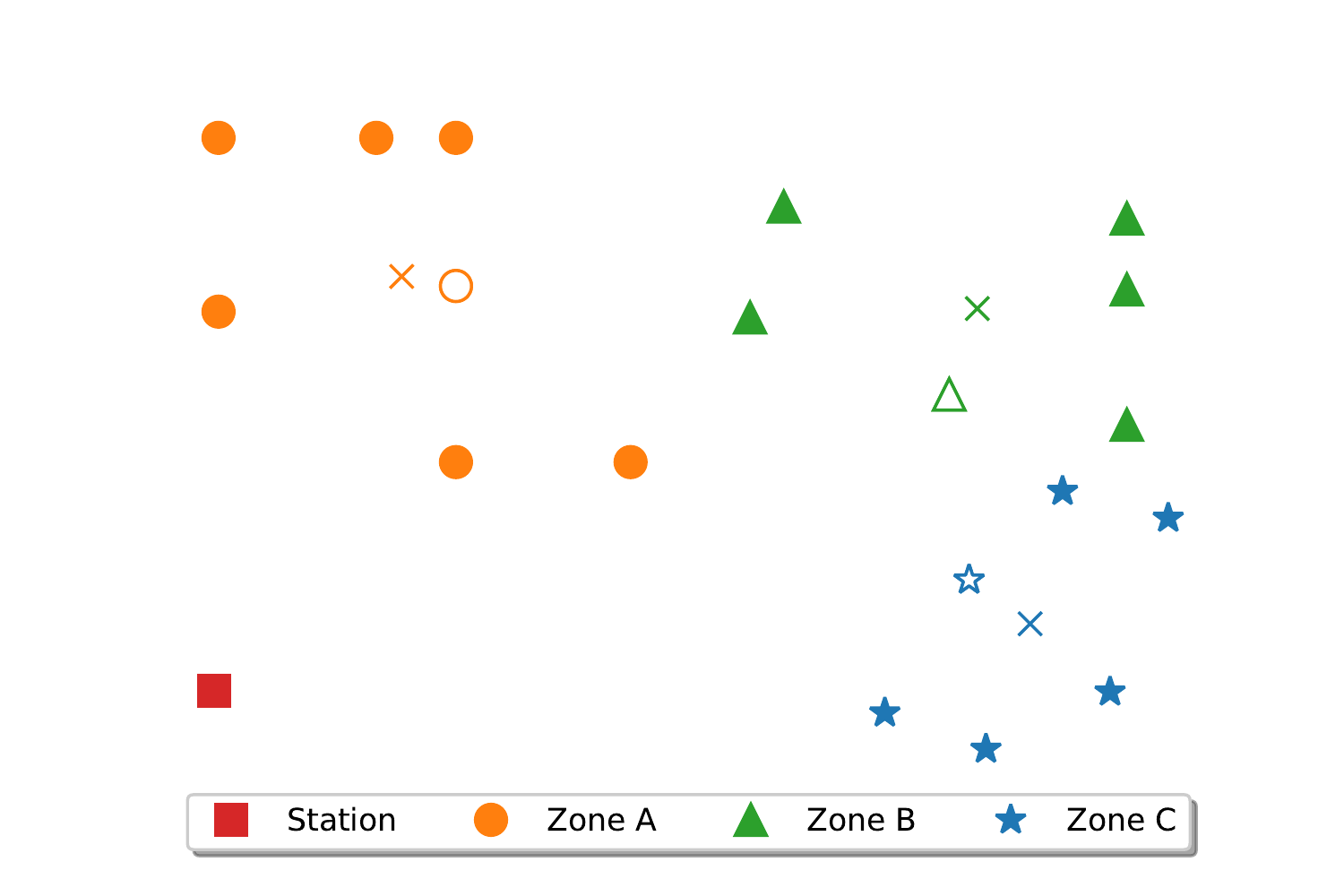}
\caption{Overview of stops with zone centers $(\times)$.}
\label{fig:method1}
\end{subfigure}
\begin{subfigure}{.49\textwidth}
\centering
\includegraphics[width=\textwidth]{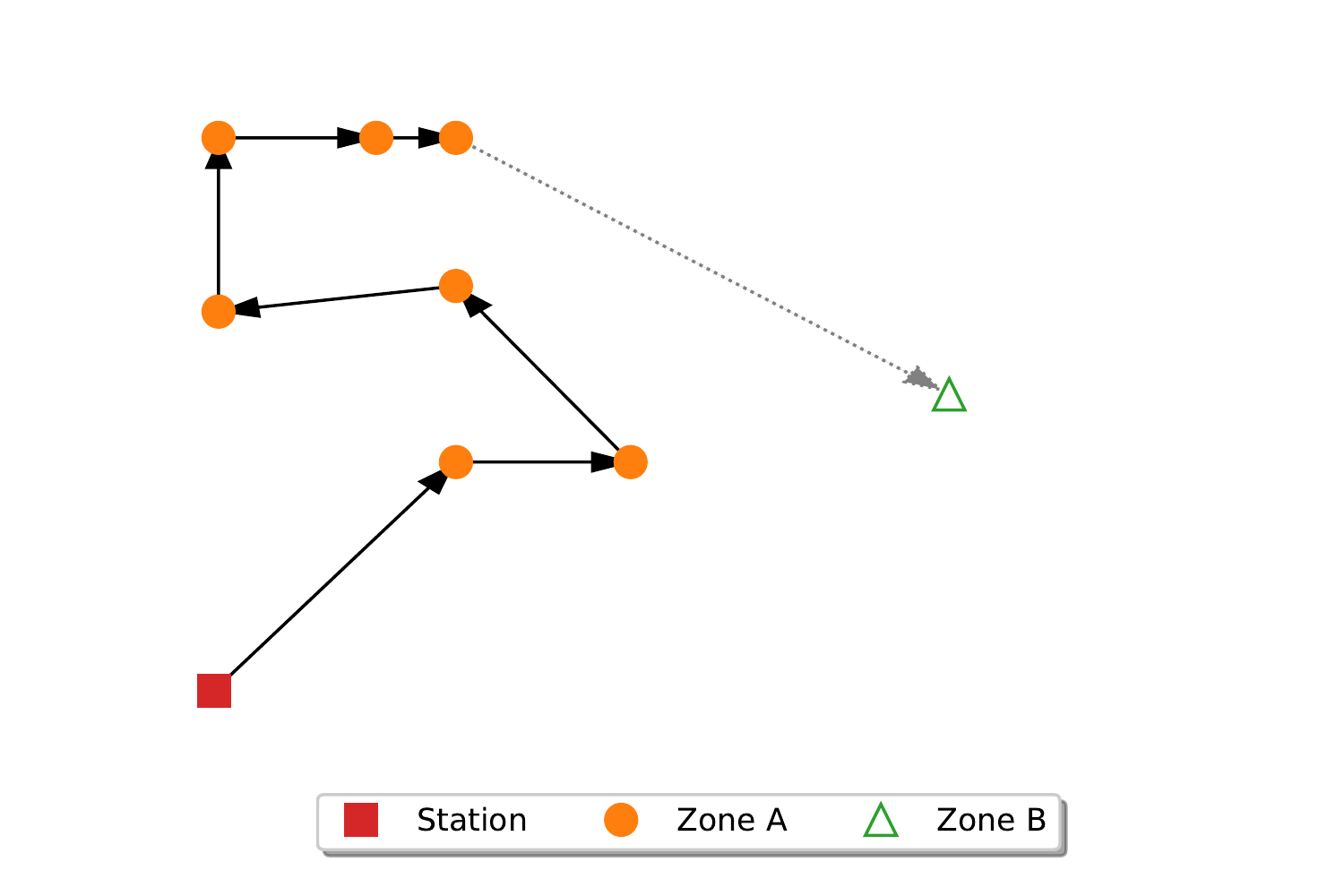}
\caption{Solution of the OTSP for zone A.}
\label{fig:method2}
\end{subfigure}
\begin{subfigure}{.49\textwidth}
\centering
\includegraphics[width=\textwidth]{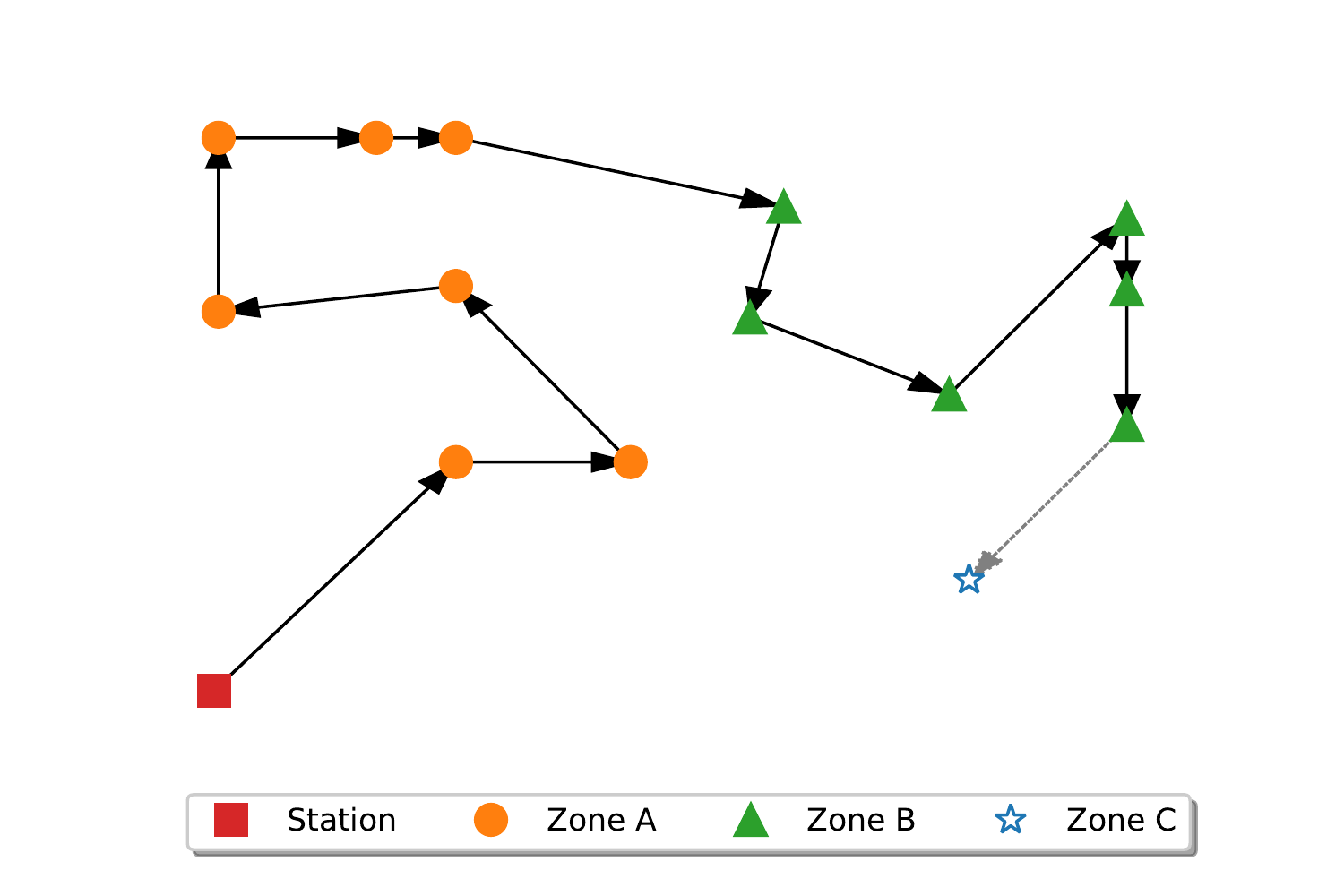}
\caption{Solution of the OTSP for zone B.}
\label{fig:method3}
\end{subfigure}
\begin{subfigure}{.49\textwidth}
\centering
\includegraphics[width=\textwidth]{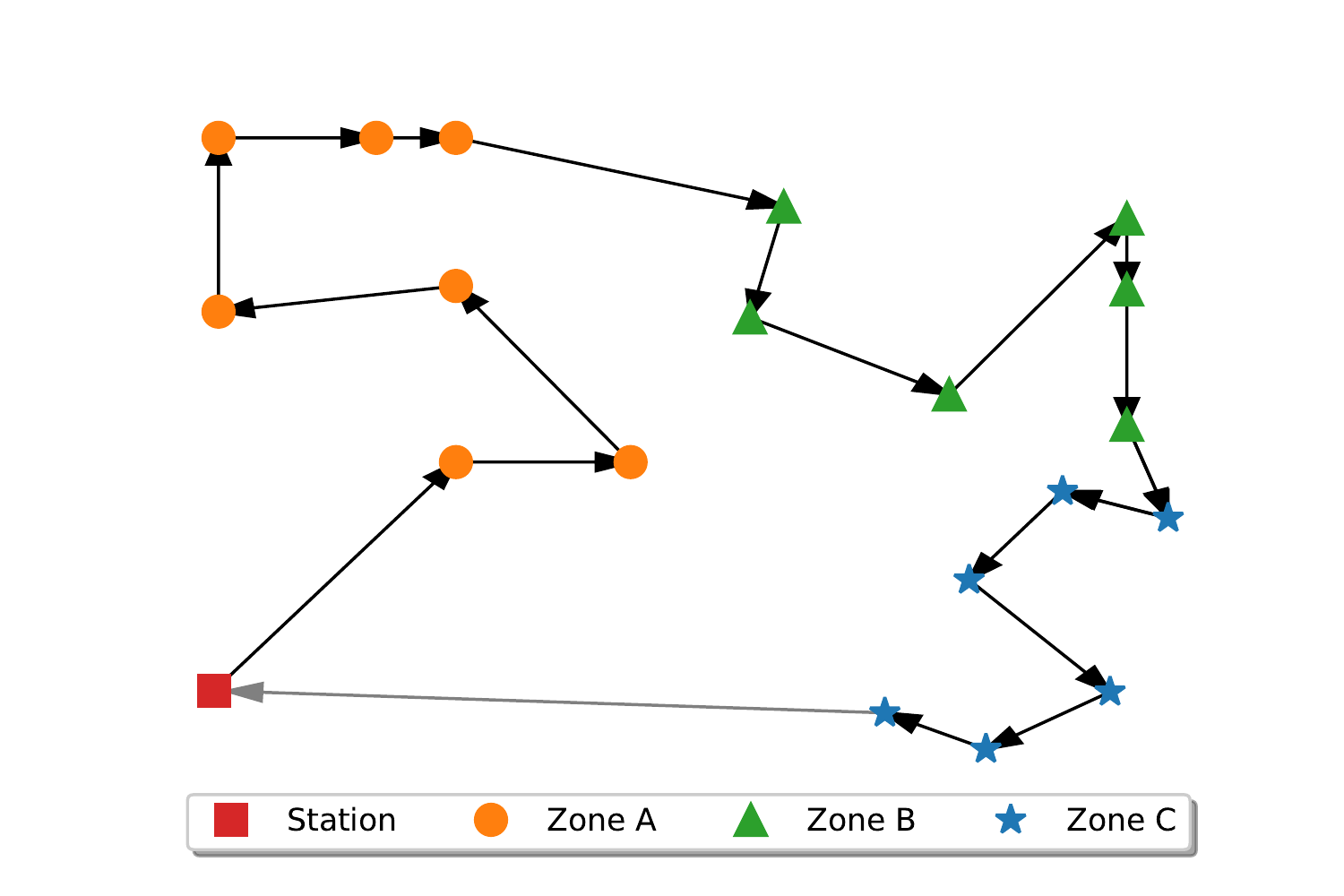}
\caption{Solution of the OTSP for zone C.}
\label{fig:method4}
\end{subfigure}
\caption{An example of the iterative procedure that is applied on (a), which each time ((b), (c) \& (d)) solves an open traveling salesman problem locally.}
\label{fig:method}
\end{figure}

At this point, we have a sequence of zones that yet has to be transformed into a sequence of stops, e.g., we are in Figure~\ref{fig:method1} where the zone sequence is $(Z_A,Z_B,Z_C).$ As the problem size in each zone is considerably smaller, drivers are able to make (near-)optimal navigational decisions that minimize travel time or distance, see~Section~\ref{sec:human_nav}. To mimic that behavior, we employ again the framework for solving TSPs, but with the additional flexibility to be an open TSP (OTSP) as we do not need to close the loop within a zone, but rather traverse from one zone to another. Although there are several options available, the procedure detailed in Figure~\ref{fig:method} is in accordance with practice. The approach builds on two elements, which are applied iteratively in each zone: 
 \begin{itemize}
     \item The last stop of a previous zone (or station, Figure~\ref{fig:method2}) is used as the starting point of the OTSP for the current zone.
     \item An additional stop of the next zone is added to each OTSP to add sense of direction to find a good stop to finish its zone. 
 \end{itemize}
As seen in Figure~\ref{fig:method1}, the extra stop (not filled) is selected as the one closest (in terms of Euclidean distance) to the computed zone centers $(\times)$. Naturally, after the solution for an OTSP is found with the additional stop added, the last stop is removed from the solution such that for the subsequent zone the starting point is again a stop in its previous zone. The procedure continues until the last zone is reached for which the station is used as the additional `stop'.

Because of clarity of the exposition, we considered Euclidean distances in the figure, but the implementation can also be fed with travel times. The resulting solution thus forms a logical sequence from a navigational point of view as it avoids path crossings and mimics human behavior when solving (O)TSPs, see Section~\ref{sec:human_nav}. A full implementation of the approach is given in the Appendix and is tailored such that it fits the dataset corresponding to the Amazon Last Mile Research Challenge. Also, we use this dataset to show the adequacy of the approach, as the achieved performance puts the approach among the top-tier of submissions for this challenge.

\section{Performance Evaluation}\label{sec:perf}
To assess the performance of the proposed methodology we employ it on a routing dataset, which was part of the Amazon Last Mile Research Challenge~\citep{routingchallenge}. Also the metrics that we use are adopted from this challenge and are provided in Section~\ref{perf}. 

The dataset counts 6,112 routes from $17$ different stations. We split the data such that 5,112 routes are used for learning and $1,000$ routes are reserved as an out-of-sample test set. Furthermore, we impute missing zone ids by taking over the id of the closest stop in terms of travel time---a sensible choice from the driver's perspective. The experiments related to our methodology are run using an Intel Core i7-10850H 2.70 GHz CPU and 32 GB RAM. The learning takes approximately 1 minute and with a test set of $1,000$ routes, it takes on average 30 minutes to render all routes. However, for choosing a right value of~$\omega$ for Eq.~\eqref{eq2} and reconsidering it in~\ref{sec:station}, we rely on 5-fold cross validation over the training set of $5,112$ routes.

\subsection{Route Prediction Metrics}\label{perf}
The route prediction performance can be evaluated in two dimensions, that is, how {\it often} and by how {\it much} a prediction $B$ deviates from the realized stop sequence $A:=(0,1,\ldots,n)$, and thus $B$ is in essence a permutation of $A$. Below, we describe these different components and how they together form an overall route performance metric. 
\begin{itemize}
    \item The {\it sequence deviation} \textbf{SD$_\text{stop}(A,B)$} measures the difference in  stop sequences $A$ and $B.$ Given that in both cases all $n+1$ stops have been visited, we take the stop sequence of $B$, and for each position we trace back when it has been visited in $A$ to create a vector $(a_0,a_1\ldots, a_n)=\pi(A)$ so that 
    $$
    \textbf{SD$_\text{stop}(A,B)$}:=\frac{2}{n(n-1)}\sum_{i=1}^n \left(|a_i-a_{i-1}|-1\right).
    $$ 
    Note that indeed if $A\equiv\pi(A)=B$, the deviation becomes $0$. Lastly, due to the fact that we deduce the zone sequence from a route by running the procedure outlined in Section~\ref{build}, we can also compute the zone sequence deviation $\textbf{SD$_\text{zone}(A,B)$}$ that on the higher level compares sequence deviations between two zone sequences $A$ and $B$.
    \item Besides considering the position in the sequence, one can also count the number of edit operations needed to come from the predicted route to the realized route, familiarly known as the \textit{Levenshtein} distance as introduced by \cite{levenshtein1966binary}. This distance metric, renamed to~\textbf{ERP$_\text{edit}(A,B)$}, counts the number of deletions and insertions to get the same sequence.
    \item Apart from considering the concurrence and concordance of stop sequences, the {\it size of deviation} between two sequences $A$ and $B$ should be evaluated as well. To this end, normalized travel time is introduced to evaluate the difference between two routes on a stop basis; let $A[i]$ and $B[j]$ be the $i$-th stop of $A$ and the $j$-th stop of $B$, with $i,j = 0,1,\dots,n$, then it is defined as
    \begin{align*}
    \textbf{time$_\text{norm}$}(A[i],B[j]) &:= Y_{A[i],B[j]} - \min_{i,j} Y_{A[i],B[j]}, \\
    Y_{A[i],B[j]} &:= \frac{t_{A[i],B[j]} - \bar{t}}{\text{std}(t)},
    \end{align*}
    where $t_{A[i],B[j]}$ is the travel time from stop $A[i]$ to stop $B[j]$, and $\bar{t}$ and $\text{std}(t)$ are the mean and standard deviation over all travel times between the stops in the route. Considering per stop the difference between two sequences, we compute the so-called Edit Distance with Real Penalty component (\textbf{ERP}$_{\textbf{norm}}$). Formally, it is computed as
    $$
    \textbf{ERP$_\text{norm}(A,B)$} := \sum_{i=0}^{n}\textbf{time$_\text{norm}$}(A[i],B[i]).
    $$
    As an aside,  dividing \textbf{ERP$_\text{norm}(A,B)$} by \textbf{ERP$_\text{edit}(A,B)$} translates to the average additional travel time incurred by deviating, and is called \textbf{ERP$_\text{ratio}(A,B)$}. 

\end{itemize}
Besides these separate metrics, they are combined to create a comprehensive score
$$
\textbf{route\_score}=\frac{\textbf{SD$_\text{stop}(A,B)$}\cdot \textbf{ERP$_\text{norm}(A,B)$}}{\textbf{ERP$_\text{edit}(A,B)$}},
$$
so that a performance score is obtained by averaging over all $I$ routes to be predicted:
\begin{align}\label{eq:perf}
\textbf{Performance}=\frac{1}{I}\sum_{i=1}^{I} \textbf{route\_score}_i.
\end{align}
Evidently, a lower route score means that the predicted sequence coincides more with the realized sequence; $0$ means they are exactly the same.

\subsection{Choosing the Weight}
\begin{figure}
\begin{tikzpicture}
\begin{axis}[legend style={legend pos={north west}},
scaled y ticks = false,
  width=\textwidth,
  height=0.5\textwidth,
legend cell align={left},
legend columns=1,
tick align=outside,
tick pos=left,
x grid style={black},
xlabel={$\omega$},
xmin=-0.01, xmax=1.01,
xtick style={color=black},
y grid style={black},
ylabel={$\textbf{Performance}$},
ymin=0.03, ymax=0.06,
ytick={0.03,0.035,0.04, 0.045,0.05,0.055,0.06},
yticklabels={0.03, ,0.04, ,0.05, ,0.06}
]

\addplot [black, thick, densely dashed, mark indices = {1, 11, 21, 31, 41, 51, 61, 71, 81, 91, 100}, mark options={solid, fill=gray}]
table {%
0 0.0528306083294636
0.01 0.0394599411073126
0.02 0.0386613227690362
0.03 0.0380899553176494
0.04 0.0377769004696848
0.05 0.0374607100351558
0.06 0.0372649634660768
0.07 0.0371087425923152
0.08 0.0370255660921876
0.09 0.036994447825233
0.1 0.0369335721909128
0.11 0.0368990908570456
0.12 0.0369092051924308
0.13 0.0368790301214092
0.14 0.036880896636853
0.15 0.0368198842273652
0.16 0.0368584358061668
0.17 0.0368164835208116
0.18 0.036853124317674
0.19 0.0367754934914966
0.2 0.0368003468909508
0.21 0.0367606154171702
0.22 0.0367892569635126
0.23 0.0367141422133666
0.24 0.036705604213122
0.25 0.0367687627386234
0.26 0.0367393044448838
0.27 0.036677948951397
0.28 0.0366342043193216
0.29 0.036655272376229
0.3 0.0366279143382564
0.31 0.03656535501935
0.32 0.0365864974415484
0.33 0.0366363862927348
0.34 0.0366005089066936
0.35 0.036632083449297
0.36 0.03658010115296
0.37 0.0365065643519006
0.38 0.0365083281255052
0.39 0.036482344102952
0.4 0.0365014228447112
0.41 0.0364439497293356
0.42 0.036478336605007
0.43 0.0365198919014704
0.44 0.036462335959333
0.45 0.0364571741162806
0.46 0.0364603569031262
0.47 0.0363839649235
0.48 0.036443772318445
0.49 0.0364339759734584
0.5 0.036423989522621
0.51 0.036336502986454
0.52 0.03633325232644
0.53 0.0363138400021084
0.54 0.0363314861698584
0.55 0.036286557251944
0.56 0.03624972232173
0.57 0.0362057819985178
0.58 0.0361928950000716
0.59 0.0361861594448964
0.6 0.0362195724364134
0.61 0.0361076525687928
0.62 0.0360970485790646
0.63 0.0360038983487932
0.64 0.035948362352335
0.65 0.0358915406194018
0.66 0.0358974639589326
0.67 0.0358636811381872
0.68 0.0358347338559262
0.69 0.0357838011052314
0.7 0.0358414280350512
0.71 0.0357361010267068
0.72 0.035742318302009
0.73 0.0356273304361988
0.74 0.0356434632920722
0.75 0.0355409732046658
0.76 0.0354375077808972
0.77 0.0352728381281284
0.78 0.0352678517970992
0.79 0.0352185893517668
0.8 0.0350821891117328
0.81 0.0350418123592986
0.82 0.0349261196119982
0.83 0.0349801573136836
0.84 0.0349012692163518
0.85 0.034824879540335
0.86 0.034703742857822
0.87 0.0346061934572966
0.88 0.034450448886237
0.89 0.0344039405889216
0.9 0.0343341690643174
0.91 0.034326552415068
0.92 0.03451345675851
0.93 0.034475834857358
0.94 0.0347409758952458
0.95 0.035494989867836
0.96 0.036716283713491
0.97 0.0383321683025456
0.98 0.0416727022652852
0.99 0.047078661232307
1 0.05502121554446
};
\addlegendentry{{\scriptsize $\omega D+ (1-\omega) (1-P)$}}

\addplot [black, thick, mark indices = {1, 11, 21, 31, 41, 51, 61, 71, 81, 91, 100}, mark options={solid,  fill=gray}]
table {%
0 0.0528274809155484
0.01 0.0391225268386316
0.02 0.0380295764635062
0.03 0.0375190129834248
0.04 0.0371550906018918
0.05 0.0370488910690808
0.06 0.0369662193807964
0.07 0.0369085382983292
0.08 0.0368283290695266
0.09 0.0367934870111216
0.1 0.036738495173183
0.11 0.036739476662057
0.12 0.0367004282035608
0.13 0.0367565495624028
0.14 0.036782039583856
0.15 0.0367805274722612
0.16 0.0366705924456536
0.17 0.0366233041510552
0.18 0.0366369306291818
0.19 0.0366819676879926
0.2 0.0366482489782432
0.21 0.036647281692223
0.22 0.0366201152689596
0.23 0.0365146172175802
0.24 0.0365631896930412
0.25 0.036516722745212
0.26 0.036487140783703
0.27 0.03650079843518
0.28 0.036582749613228
0.29 0.036506086619738
0.3 0.0365061017483352
0.31 0.0364865017109794
0.32 0.036502086698499
0.33 0.0364443595214868
0.34 0.0364136981917818
0.35 0.0364326944013894
0.36 0.036430184422508
0.37 0.0363711375631908
0.38 0.0363149591849294
0.39 0.0363472568694692
0.4 0.036332154973839
0.41 0.0362898887952754
0.42 0.0362768560454236
0.43 0.036196529508128
0.44 0.036238504991336
0.45 0.0361595404800158
0.46 0.0361763962524094
0.47 0.036167017804672
0.48 0.036118378623624
0.49 0.0361039303401852
0.5 0.0359735650206936
0.51 0.0359902645603182
0.52 0.0360172546694978
0.53 0.0359179761857564
0.54 0.0359884261979966
0.55 0.0359229391859504
0.56 0.0359415600955062
0.57 0.0358969664276412
0.58 0.03588054798185
0.59 0.035824231281616
0.6 0.0357908015333894
0.61 0.0358035490407964
0.62 0.0357021535623112
0.63 0.0356261364820556
0.64 0.0355698863261022
0.65 0.0354549320001984
0.66 0.0353297343315866
0.67 0.0352520312179492
0.68 0.035169537702474
0.69 0.035110587443023
0.7 0.035058502087321
0.71 0.0349284806710404
0.72 0.0347676895118138
0.73 0.034605683534399
0.74 0.034564538743834
0.75 0.0344735014108372
0.76 0.0343784249771266
0.77 0.034359355106931
0.78 0.0341827551953038
0.79 0.034154401322084
0.8 0.0339069848023312
0.81 0.0338127538297612
0.82 0.03375194235606
0.83 0.0336041742726576
0.84 0.0334533029750376
0.85 0.0332366518372262
0.86 0.0331462316571892
0.87 0.0328671383744084
0.88 0.0327931564546308
0.89 0.0328096407640982
0.9 0.0326647650939902
0.91 0.0326467387522282
0.92 0.0329940850785976
0.93 0.0331357174247622
0.94 0.0336579678155406
0.95 0.0344319556824146
0.96 0.0355385143176186
0.97 0.037413868723277
0.98 0.0409380363402546
0.99 0.0467336615686346
1 0.0529658982151836
};
\addlegendentry{{\scriptsize $\omega T+ (1-\omega) (1-P)$}}
\end{axis}
\end{tikzpicture}
  \caption{Comparison of performances, as defined in Eq~\eqref{eq:perf}, when varying $\omega$ in Eq.~\eqref{eq2}, where Euclidean distance is dashed and (expected) travel time solid. The curves are obtained averaging over the 5-fold cross-validation scores.}
  \label{fig:structure}
\end{figure}
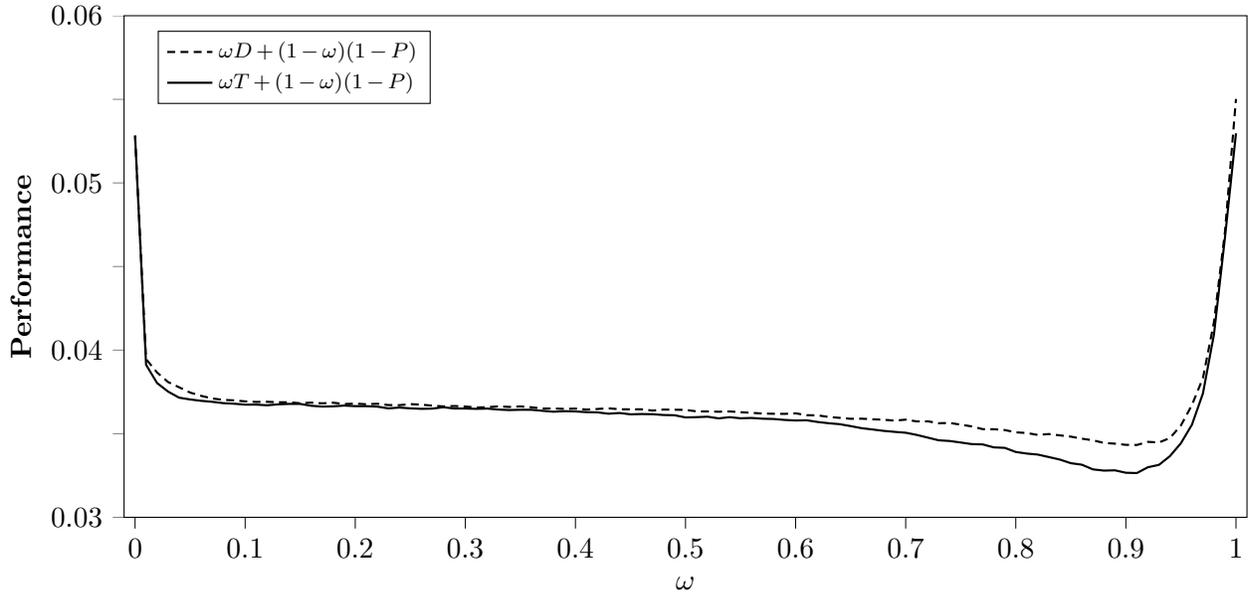

As a first experiment, we consider varying the weight parameter $\omega$ in Eq.~\eqref{eq2}. Furthermore, we can replace distance by travel times $T$ as these are provided in Amazon's dataset. The performances when varying $\omega$ in both versions are displayed in Figure~\ref{fig:structure}.

As expected, we find that only relying on historical information ($\omega=0$), or only using distance or travel times ($\omega=1$) does not perform as well as having a trade-off between the two. Overall, we observe that an $\omega$ value around $0.9$ achieves the best trade-off, but actually any value between 0.1 and 0.9 in~Eq.~\eqref{eq2} (dashed and solid lines) will result in a performance of around $0.035.$ This observation demonstrates the robustness of the model as being insensitive to misspecification of $\omega$. In addition, we conclude that the difference between travel times or Euclidean distance is small and thus distance can be used as a good proxy for travel time.

\subsection{Adjusting the Station's Weights}\label{sec:station}
\begin{figure}
    \centering
    \includegraphics[width=.75\textwidth]{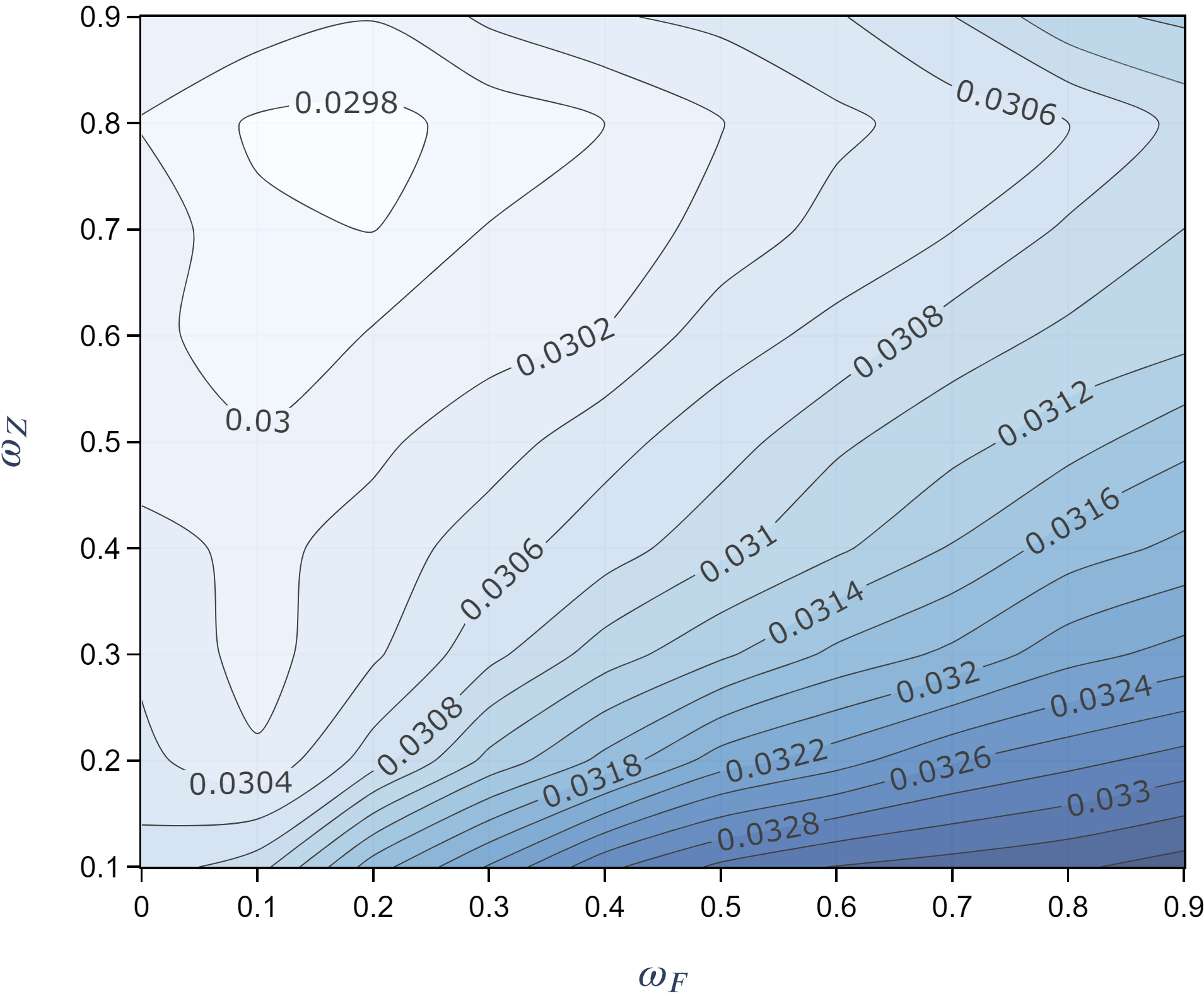}
    \caption{A contour plot showing the performance across different weightings of $\omega_F$ (station-to-zone) and $\omega_Z$ (zone-to-zone) transitions, while keeping $\omega_L=1$ (zone-to-station) in Eq.~\eqref{eq3}. The values are obtained averaging over the 5-fold cross-validation scores.}
    \label{fig:contour}
\end{figure}

Reconsidering Eq.~\eqref{eq2}, wherein we define an $\Omega$ matrix to weigh historical information against travel times, we further study the role of the station as the starting and ending point of a route. For the sake of readability, we consider a matrix $\Omega^{(r)}$ which contains the weights in $\Omega$ corresponding to the station and zones that are visited in a specific route $r$ where the station is related to the first row and column. So as to refine the role of the station, we define $\Omega^{(r)}$ matrix as
\[
\Omega^{(r)} = 
\left(
\begin{array}{c|cccc}
0  & \omega_F & \dots & \omega_F \\ \hline
\omega_L & \omega_Z & \dots & \omega_Z \\
\vdots & \vdots & \ddots & \vdots \\
\omega_L & \omega_Z & \dots & \omega_Z
\end{array}
\right),
\]
where $\omega_F \in [0,1]$ is the weight in the cost function from the station to the \textit{first} zone. Analogously, $\omega_Z \in [0,1]$ sets the balance for \textit{zone-to-zone} transitions, and \textit{lastly} $\omega_L \in [0,1]$ is the weight corresponding to the return to the station. The elements $\Omega_{ij}$ from this matrix are put in the following refinement of~Eq.~\eqref{eq2} with travel times: 
\begin{equation}\label{eq3}
C_{ij} = \Omega_{ij}\, T_{ij} + (1 - \Omega_{ij})\,(1 - P_{ij}),
\end{equation}
where, similarly to Eq.~\eqref{eq2}, \tc{the diagonal values are set to zero. With this generalized cost matrix we first set $\omega_L=1$ as we anticipate that at the end of a route heading back to the station should be the only concern, and thus the number of times that a specific zone acted as the finish point of a route can be disregarded. This resonates with the thought that a driver is only focused at the end on the travel time it takes to return to station. Having this parameter restrained, we study combinations of $\omega_F$ and $\omega_Z$ to obtain the contour plot in Figure~\ref{fig:contour}. In the plot, we observe that the optimal point should lie somewhere in the region where $\omega_F\in[0.1,0.2]$ and $\omega_Z\in[0.7,0.8]$. So when starting a route in the first transition from the station, historical information should be weighted more important than regular zone-to-zone transitions. This implicates that when determining the route globally preference and navigational considerations have a greater impact than the standard trade-off with travel times.} 


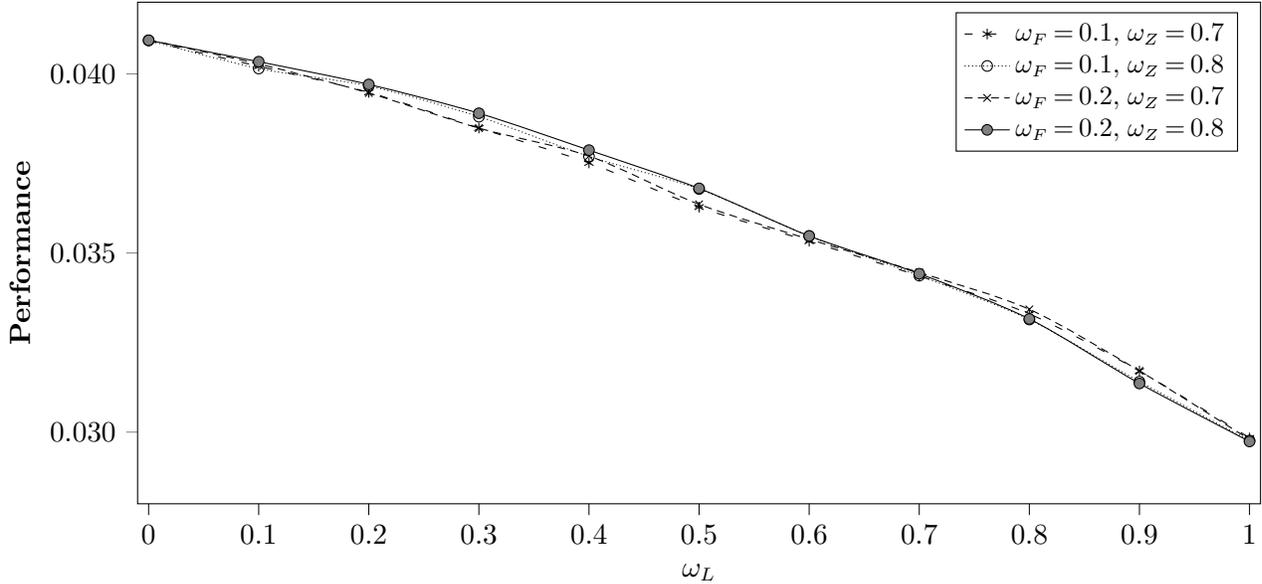
\begin{figure}
    \centering
      \begin{tikzpicture}
    \begin{axis}[legend style={font=\small},
    legend cell align={left},
scaled y ticks = false,
  width=\textwidth,
  height=0.5\textwidth,
legend cell align={left},
tick align=outside,
tick pos=left,
x grid style={black},
xlabel={$\omega_L$},
xmin=-0.01, xmax=1.01,
xtick style={color=black},
y grid style={black},
ylabel={$\textbf{Performance}$},
ymin=0.028, ymax=0.042,
ytick={0.03,0.035,0.04},
yticklabels={0.030,0.035,0.040},
    cycle list name=black white,
    smooth]

\addlegendentry{$\omega_F = 0.1,\, \omega_Z = 0.7$}  
     \addplot+[color=black,
                style={loosely dashed},
                mark=asterisk,
                mark options={solid, fill=gray}] 
table {%
0 0.0409198675129072
0.1 0.040212320835062
0.2 0.0394825002714804
0.3 0.0384883739370356
0.4 0.0375133092416708
0.5 0.0362804724020588
0.6 0.0353273113137516
0.7 0.0343578374156526
0.8 0.0332996705731322
0.9 0.0317061071700078
1 0.0298383862045326
};

\addlegendentry{$\omega_F = 0.1,\, \omega_Z = 0.8$}  
         \addplot+[color=black,
                style={densely dotted},
                mark=o,
                mark options={solid, fill=gray}]
table {%
0 0.040932836618523
0.1 0.0401483100612138
0.2 0.039660489608073
0.3 0.0388048040143306
0.4 0.0376964398837066
0.5 0.0367860153461906
0.6 0.0354710469196186
0.7 0.0343677153702544
0.8 0.033142420273566
0.9 0.0314089420225934
1 0.0297662620253192
};

\addlegendentry{$\omega_F = 0.2,\, \omega_Z = 0.7$}
     \addplot+[color=black,
                style={densely dashed},
                mark=x,
                mark options={solid, fill=gray}]
table {%
0 0.0409445256688752
0.1 0.0402766829481996
0.2 0.0394649021974672
0.3 0.038485468533634
0.4 0.0377055968254196
0.5 0.0363592037011954
0.6 0.035383327273823
0.7 0.0344560132737094
0.8 0.0334267509825942
0.9 0.0317087685670362
1 0.0297933749515538
};

\addlegendentry{$\omega_F = 0.2,\, \omega_Z = 0.8$}
     \addplot+[color=black,
   style=solid,
                mark=*,
                mark options={solid, fill=gray}]
table {%
0 0.0409370182460814
0.1 0.0403392089888922
0.2 0.0397044589710658
0.3 0.0389004502091348
0.4 0.0378707519712416
0.5 0.0368034546934168
0.6 0.035476567047145
0.7 0.034421504293425
0.8 0.0331599107219942
0.9 0.031358367749612
1 0.0297396275795116
};
\end{axis}
    \end{tikzpicture}
    \caption{A sensitivity study on $\omega_L$, (zone-to-station), for selected combinations of $\omega_F$ and $\omega_Z$ that lie around the inner contour of~Figure~\ref{fig:contour}. The performances, as defined in Eq.~(\ref{eq:perf}), are obtained by averaging over the 5-fold cross-validation scores.}
    \label{fig:omegaL}
\end{figure}

The choice of $\omega_L=1$ is justified by the degradation in performance for values less than 1. This is shown in Figure~\ref{fig:omegaL} by taking four combinations of $(\omega_F,\omega_Z)$ that lie around the inner contour---close to the optimum. 
In all, when comparing the points in detail, we conclude that the best performance in our cross-validation is obtained when $(\omega_F,\omega_Z,\omega_L)=(0.2,0.8,1)$, which results in the reported performance of the final model in~Figure~\ref{fig:route_score}. 

\subsection{Comparison against Benchmarks}
\tc{In Table~\ref{tab:baseline}, we provide the performance of our approach abbreviated to \textit{LG-OL}  and compare it to two benchmarks to shed some light onto the range of performance scores that can be attained. Besides the metrics outlined in Section~\ref{sec:perf} we also added the mean travel time in seconds in the  \textbf{Time} column. First, we apply the Nearest Neighbor algorithm (Figure~\ref{fig:nn_solution}) on the out-of-sample test set of $1,000$ routes. This algorithm sequentially adds a connection from the last stop to the nearest unvisited stop; a strategy that is sometimes considered to model human route planning behavior~\citep{graham2000traveling,wiener2009planning}. As another idea, consider the route which minimizes the travel times as the route a driver follows (Figure~\ref{fig:tsp_solution}). Such route is generally hard to compute, but within the range of state-of-the-art solvers (e.g., Gurobi). Surprising is that drivers average travel time (9496 seconds) is about 18\% more than when travel time minimization was the only objective (Full TSP); indeed showing that drivers' navigational decisions cannot be approached as a standard optimization problem.}

\begin{figure}
    \centering

\begin{tikzpicture}

\definecolor{color0}{rgb}{0.12156862745098,0.466666666666667,0.705882352941177}
\definecolor{color1}{rgb}{1,0.498039215686275,0.0549019607843137}
\definecolor{color2}{rgb}{0.172549019607843,0.627450980392157,0.172549019607843}

\begin{groupplot}[group style={group size=2 by 1}]
\nextgroupplot[
tick align=outside,
tick pos=left,
x grid style={white!69.0196078431373!black},
xmin=0.5, xmax=3.5,
xtick style={color=black},
xtick={1,2,3},
x post scale=2,
xticklabels={Nearest Neighbor, Full TSP, LG-OL},
y grid style={white!69.0196078431373!black},
ylabel={Route Score},
ymin=-0.0117803820966705, ymax=0.256577275978241,
ytick style={color=black},
ytick={0,0.05,0.1,0.15,0.2,0.25},
yticklabels={0.00,0.05,0.10,0.15,0.20,0.25}
]
\addplot [black]
table {%
0.85 0.068867325063026
1.15 0.068867325063026
1.15 0.139999565325417
0.85 0.139999565325417
0.85 0.068867325063026
};
\addplot [black]
table {%
1 0.068867325063026
1 0.010150614746743
};
\addplot [black]
table {%
1 0.139999565325417
1 0.2443792006112
};
\addplot [black]
table {%
0.925 0.010150614746743
1.075 0.010150614746743
};
\addplot [black]
table {%
0.925 0.2443792006112
1.075 0.2443792006112
};
\addplot [black]
table {%
1.85 0.0472752191098515
2.15 0.0472752191098515
2.15 0.104285198970099
1.85 0.104285198970099
1.85 0.0472752191098515
};
\addplot [black]
table {%
2 0.0472752191098515
2 0.005260057624565
};
\addplot [black]
table {%
2 0.104285198970099
2 0.187429004095569
};
\addplot [black]
table {%
1.925 0.005260057624565
2.075 0.005260057624565
};
\addplot [black]
table {%
1.925 0.187429004095569
2.075 0.187429004095569
};
\addplot [black]
table {%
2.85 0.00581107555514675
3.15 0.00581107555514675
3.15 0.0422827657254975
2.85 0.0422827657254975
2.85 0.00581107555514675
};
\addplot [black]
table {%
3 0.00581107555514675
3 0.000417693270371
};
\addplot [black]
table {%
3 0.0422827657254975
3 0.094856919500038
};
\addplot [black]
table {%
2.925 0.000417693270371
3.075 0.000417693270371
};
\addplot [black]
table {%
2.925 0.094856919500038
3.075 0.094856919500038
};
\addplot [black]
table {%
0.85 0.09853286317581
1.15 0.09853286317581
};
\addplot [black]
table {%
1.85 0.071217842337108
2.15 0.071217842337108
};
\addplot [black]
table {%
2.85 0.015468360116014
3.15 0.015468360116014
};
\end{groupplot}
\end{tikzpicture}
    \caption{Boxplot representing the route score distribution on the out-of-sample test set of 1,000 routes for our model, and the two benchmarks of Table~\ref{tab:baseline}.
    \label{fig:route_score}}
\end{figure}
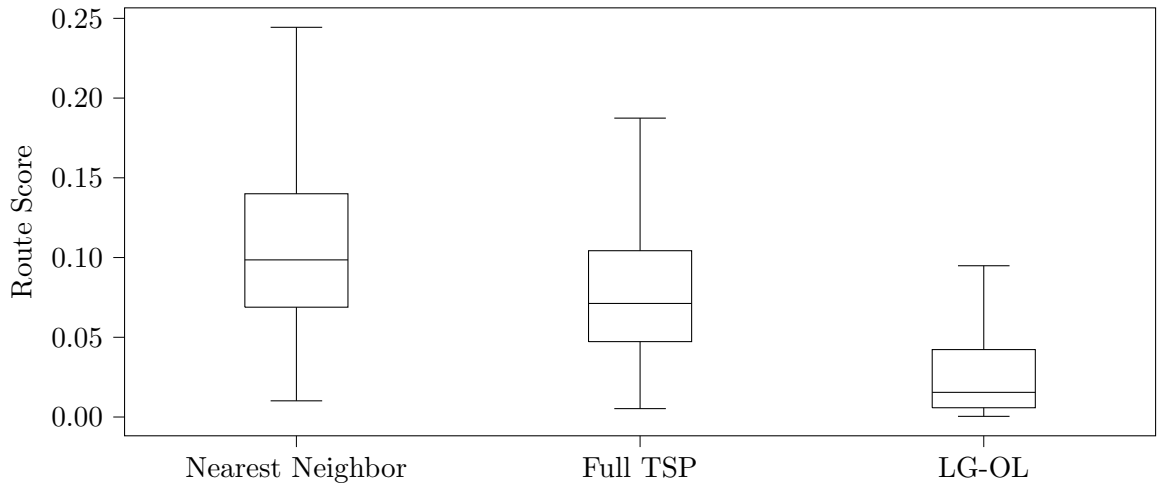

\tc{Our methodology obtains a prediction performance (\textbf{Performance}) of $0.0299$---a more than $64\%$ reduction compared to the aforementioned benchmarks. Interestingly, the Nearest Neighbor and full TSP models have similar \textbf{ERP$_\text{ratio}$} scores, but the full TSP renders better sequences (lower \textbf{SD$_\text{zone}$} and \textbf{SD$_\text{stop}$} scores). When considering our individual route scores of the test set in Figure~\ref{fig:route_score} against the benchmarks, we find that our model succeeds to have $40\%$ of the test routes to be near-perfectly predicted, as they are below 0.01; and even gets $80\%$ of them below the $0.05$ threshold. In retrospect, reconsidering Figure~\ref{fig:structure} we find that even for any $\omega$ the unadjusted approach outperforms these benchmarks, supporting the fundamental idea of our learn global, optimize local approach.}

\tc{Evaluating our performance to the submissions for the challenge~\citep{routingchallenge2}, we find that the approach propels itself in the top-tier of submissions, which performances range from 0.025  to 0.050. Finally, to verify the correctness and potential of the approach, we  did an additional hypothetical experiment on the test set of which the results are provided in Table~\ref{tab:baseline} as \textit{LG-OL (hypothetical)}. We extracted the zone sequences using Section~\ref{build} on the test set so as to consider the situation that we are fully capable to predict the zone sequence right. Thus remaining are the inconsistencies between the design choices made in Section~\ref{sec:model},  and drivers' behavior that is not reflected in the methodology. In the experiment, we find that the performance becomes 0.0094, which is well below all submissions for this challenge.} 



\begin{table}
\TABLE{Comparison of the methodology against benchmarks on the out-of-sample test set.}
{\begin{tabular}{@{}lrrrrr@{}}
\toprule
\textbf{Model}  & \textbf{SD$_\text{zone}$} & \textbf{SD$_\text{stop}$} & \textbf{ERP$_\text{ratio}$} & \tc{\textbf{Time} (s)} & \textbf{Performance} \\ \midrule
Nearest Neighbor & 0.2441 & 0.0873 & 1.3382 & 9865 & 0.1119 \\
Full TSP         & 0.1912 & 0.0650 & 1.3330 & 8036 & 0.0826 \\ 
LG-OL            & 0.0940 & 0.0335 & 0.7389 & 8994 & 0.0299 \\ \hline
\textit{LG-OL (hypothetical)}           & 0 & 0.0238 & 0.3831 & 9535 & 0.0094 \\
\textit{Driver}           & 0 & 0 & 0 & 9496 & 0 \\ \bottomrule
\end{tabular}}
{\label{tab:baseline}}
\end{table}

\section{Conclusion and Discussion}\label{sec:concl}
The costs associated with last-mile delivery surmount any other shipping costs. Hence, good predictions about the route taken are crucial from an operations management point of view. The proposed approach leads to several broadly applicable insights about delivery driving. 

The approach first determines a zone sequence by solving a traveling salesman problem (TSP) over a cost matrix which is a weighted combination of historical information, which has to be learnt from realized routes, and distance (or travel time) between zones. Next, on the local level, adhering to the established zone sequence on the global level, the order of stops to be visited in a zone is found by solving an open TSP (OTSP), which generates locally optimal solutions. Repeating this for all subsequent zones and patching them together generates a feasible route. Therefore this approach can be considered to learn global (over zones) and optimize local (within zones). Besides that the approach can be argued to follow the logic of human navigation, the approach is convenient, because it only requires learning on the zone level as by tracking the zone transitions. The breakdown makes the approach relatively fast, as the zone-sequence TSP and within-zone OTSPs are not so computationally involved compared to solving a single TSP for the entire route. This also renders the approach amenable to real-time training using newly realized routes of the same day. 

The key element in the approach is the computation of a cost matrix that combines distance and historical information when establishing the zone sequence as a TSP. The experiments show that a weighted combination with travel times provides a slight performance improvement over using Euclidean distance. This fact underpins that when determining the zone sequence an implicit trade-off is made between the driver's zone-to-zone preferences and a measure of distance. Studying the weight value by changing it from $0$ to $1,$ where the setting of 1 corresponds to only considering travel time, uncovers that the exact value is not critical. Overall, the best performance is obtained when the value is close to $0.9$. However, the weighting of station-to-zone and zone-to-station transitions can be further adapted in the cost matrix to better mimic driver's behavior resulting in another performance improvement. 
In fact, these adaptations are motivated by the logic that at the beginning of a route the driver is more concerned with where to start, whereas at the end, when all stops are visited, the only consideration is the return to the station. Indeed, the experiments confirm these hypotheses and we find that the travel time matrix  should be weighted less important in the station-to-zone transitions than in the zone-to-zone transitions, and from zone-to-station it is the only factor. \tc{These insights come with the caveat that drivers were suggested a route beforehand~\citep[see][]{routingchallenge}, which might have tainted the driver's decisions.}  

\tc{The approach as-is only utilizes historical stop data (with zone ids) in making predictions. Using additional data which are part of the case study, such as package size, start time, quality of the route, did not provide any clear starting points for improvement. Since the zone ids were given, it is worthwhile to further study the underlying segmentation, because it determines the breakdown of the problem into a global and a series of local problems and as such it is pivotal to the performance of the approach. We have seen that occurrence of zone revisitations are rare. However, in 2.6\% of the routes there are more than five revisitations, which hints that in these routes it might be valuable to  reconsider the zone segmentation.}

\tc{The test set also provides us an opportunity to quantify the remaining learning gap. By first extracting the zone sequences in the test set after which we apply the approach, we arrive at a score just below $0.01.$ This shows that the learn global and optimize local approach can yield a score as low as 0.009 (as a reference we get it to 0.0299, as shown in Section~\ref{sec:perf}). The learning took place on only 5,112 routes, which results in a sparse transition matrix as the zone transition probabilities, when available, are learnt based on the relatively few routes available for each region.
Increasing the number of routes for learning, for example by daily tracking of realized routes and updating the count matrix accordingly, the approach can likely improve its accuracy on zone sequence prediction, which have a profound impact on the overall performance. }  


Another endeavor is to further study the cost matrix structure. The methodology readily incorporates other functional forms, and also different weightings between historic information and distance-based measures. These weightings might depend on the size of the dataset available and even on the length of the route itself. Nevertheless, we find that a linear structure between distance or travel time and historical information is robust to misspecification of this weight parameter and is practically sound.

Finally, the methodology can be employed in other (last-mile) delivery concepts, see for example~\cite{boysen2021last}, or generalized to other operational optimization problems, which are susceptible to human interference in the decision making process---a logical starting point is the allied Vehicle Routing Problem, e.g.,~\cite{braysy2005vehicle}. In all, this research demonstrates that integrating learning in an optimization pipeline leverages \tc{its connection and value to practice.} 


\clearpage


%

%
%


\bibliographystyle{informs2014trsc} 
\bibliography{mybibfile} 

\begin{thebibliography}{48}
\providecommand{\natexlab}[1]{#1}
\providecommand{\url}[1]{\texttt{#1}}
\providecommand{\urlprefix}{URL }

\bibitem[{{Amazon and MIT Center for Transportation \&
  Logistics}(2021{\natexlab{a}})}]{routingchallenge}
{Amazon and MIT Center for Transportation \& Logistics}, 2021{\natexlab{a}}
  \emph{{Amazon Last-Mile Routing Research Challenge}}.
  \urlprefix\url{https://routingchallenge.mit.edu/about-the-challenge/},
  (accessed on 30 April 2022).

\bibitem[{{Amazon and MIT Center for Transportation \&
  Logistics}(2021{\natexlab{b}})}]{routingchallenge2}
{Amazon and MIT Center for Transportation \& Logistics}, 2021{\natexlab{b}}
  \emph{{Amazon Last-Mile Routing Research Challenge}}.
  \urlprefix\url{https://routingchallenge.mit.edu/last-mile-routing-challenge-team-performance-and-leaderboard/},
  (accessed on 30 April 2022).

\bibitem[{Applegate et~al.(2011)Applegate, Bixby, Chv{\'a}tal,
  \protect\BIBand{} Cook}]{applegate2011traveling}
Applegate DL, Bixby RE, Chv{\'a}tal V, Cook WJ, 2011 \emph{The Traveling
  Salesman Problem} (Princeton University Press).

\bibitem[{Boyer, Prud'homme, \protect\BIBand{} Chung(2009)}]{boyer2009last}
Boyer KK, Prud'homme AM, Chung W, 2009 \emph{The last mile challenge:
  Evaluating the effects of customer density and delivery window patterns}.
  \emph{Journal of Business Logistics} 30(1):185--201.

\bibitem[{Boysen, Fedtke, \protect\BIBand{}
  Schwerdfeger(2021)}]{boysen2021last}
Boysen N, Fedtke S, Schwerdfeger S, 2021 \emph{Last-mile delivery concepts: a
  survey from an operational research perspective}. \emph{OR Spectrum}
  43(1):1--58.

\bibitem[{Br{\"a}ysy \protect\BIBand{}
  Gendreau(2005{\natexlab{a}})}]{braysy2005vehicle}
Br{\"a}ysy O, Gendreau M, 2005{\natexlab{a}} \emph{Vehicle routing problem with
  time windows, part {I}: Route construction and local search algorithms}.
  \emph{Transportation Science} 39(1):104--118.

\bibitem[{Br{\"a}ysy \protect\BIBand{}
  Gendreau(2005{\natexlab{b}})}]{braysy2005vehicle2}
Br{\"a}ysy O, Gendreau M, 2005{\natexlab{b}} \emph{Vehicle routing problem with
  time windows, part {II}: Metaheuristics}. \emph{Transportation Science}
  39(1):119--139.

\bibitem[{Campuzano, Obreque, \protect\BIBand{}
  Aguayo(2020)}]{campuzano2020accelerating}
Campuzano G, Obreque C, Aguayo MM, 2020 \emph{Accelerating the
  {M}iller-{T}ucker-{Z}emlin model for the asymmetric traveling salesman
  problem}. \emph{Expert Systems with Applications} 148:113229.

\bibitem[{Canoy \protect\BIBand{} Guns(2019)}]{canoy2019vehicle}
Canoy R, Guns T, 2019 \emph{Vehicle routing by learning from historical
  solutions}. \emph{International Conference on Principles and Practice of
  Constraint Programming}, 54--70 (Springer).

\bibitem[{Chopra(2003)}]{chopra2003designing}
Chopra S, 2003 \emph{Designing the distribution network in a supply chain}.
  \emph{Transportation Research Part E: Logistics and Transportation Review}
  39(2):123--140.

\bibitem[{Current \protect\BIBand{} Schilling(1989)}]{CURSHI89}
Current JR, Schilling DA, 1989 \emph{The covering salesman problem}.
  \emph{Transportation Science} 23(3):208--213.

\bibitem[{Dantzig \protect\BIBand{} Ramser(1959)}]{Dantzig1959}
Dantzig GB, Ramser JH, 1959 \emph{The truck dispatching problem}.
  \emph{Management Science} 6(1):80--91.

\bibitem[{Deloison et~al.(2020)Deloison, Hannon, Huber, Heid, Klink, Sahay,
  \protect\BIBand{} Wolff}]{deloison2020future}
Deloison T, Hannon E, Huber A, Heid B, Klink C, Sahay R, Wolff C, 2020
  \emph{The future of the last-mile ecosystem: Transition roadmaps for
  public-and private-sector players} (World Economic Forum).

\bibitem[{Dry et~al.(2006)Dry, Lee, Vickers, \protect\BIBand{}
  Hughes}]{dry2006human}
Dry M, Lee MD, Vickers D, Hughes P, 2006 \emph{Human performance on visually
  presented traveling salesperson problems with varying numbers of nodes}.
  \emph{The Journal of Problem Solving} 1(1):4.

\bibitem[{eMarketer(2020)}]{eMarketer2020}
eMarketer, 2020 \emph{{US} ecommerce 2020, coronavirus boosts ecommerce
  forecast and will accelerate channel-shift}.
  \urlprefix\url{https://www.emarketer.com/content/us-ecommerce-2020},
  (accessed on 30 April 2022).

\bibitem[{Gendreau et~al.(1998)Gendreau, Hertz, Laporte, \protect\BIBand{}
  Stan}]{gendreau1998generalized}
Gendreau M, Hertz A, Laporte G, Stan M, 1998 \emph{A generalized insertion
  heuristic for the traveling salesman problem with time windows}.
  \emph{Operations Research} 46(3):330--335.

\bibitem[{Gendreau, Laporte, \protect\BIBand{}
  Vigo(1999)}]{gendreau1999heuristics}
Gendreau M, Laporte G, Vigo D, 1999 \emph{Heuristics for the traveling salesman
  problem with pickup and delivery}. \emph{Computers \& Operations Research}
  26(7):699--714.

\bibitem[{Gevaers, Van~de Voorde, \protect\BIBand{}
  Vanelslander(2011)}]{gevaers2011characteristics}
Gevaers R, Van~de Voorde E, Vanelslander T, 2011 \emph{Characteristics and
  typology of last-mile logistics from an innovation perspective in an urban
  context}. \emph{City Distribution and Urban Freight Transport} (Edward Elgar
  Publishing).

\bibitem[{Gevaers, {Van de Voorde}, \protect\BIBand{}
  Vanelslander(2014)}]{GEVAERS2014398}
Gevaers R, {Van de Voorde} E, Vanelslander T, 2014 \emph{Cost modelling and
  simulation of last-mile characteristics in an innovative b2c supply chain
  environment with implications on urban areas and cities}. \emph{Procedia -
  Social and Behavioral Sciences} 125:398--411, eighth International Conference
  on City Logistics 17-19 June 2013, Bali, Indonesia.

\bibitem[{Goodman(2005)}]{goodman2005whatever}
Goodman RW, 2005 \emph{Whatever you call it, just don't think of last-mile
  logistics, last}. \emph{Global Logistics \& Supply Chain Strategies} 9(12).

\bibitem[{Graham, Joshi, \protect\BIBand{} Pizlo(2000)}]{graham2000traveling}
Graham SM, Joshi A, Pizlo Z, 2000 \emph{The traveling salesman problem: A
  hierarchical model}. \emph{Memory \& Cognition} 28(7):1191--1204.

\bibitem[{Gunasekaran, Patel, \protect\BIBand{}
  McGaughey(2004)}]{gunasekaran2004framework}
Gunasekaran A, Patel C, McGaughey RE, 2004 \emph{A framework for supply chain
  performance measurement}. \emph{International Journal of Production
  Economics} 87(3):333--347.

\bibitem[{Jiang et~al.(2014)Jiang, Gao, Li, Wu, \protect\BIBand{}
  Pei}]{jiang2014hierarchical}
Jiang J, Gao J, Li G, Wu C, Pei Z, 2014 \emph{Hierarchical solving method for
  large scale tsp problems}. \emph{International symposium on neural networks},
  252--261 (Springer).

\bibitem[{Karp(1977)}]{karp1977probabilistic}
Karp RM, 1977 \emph{Probabilistic analysis of partitioning algorithms for the
  traveling-salesman problem in the plane}. \emph{Mathematics of Operations
  Research} 2(3):209--224.

\bibitem[{Kool, van Hoof, \protect\BIBand{} Welling(2019)}]{kool2019attention}
Kool W, van Hoof H, Welling M, 2019 \emph{Attention, learn to solve routing
  problems!} \emph{International Conference on Learning Representations},
  \urlprefix\url{https://openreview.net/forum?id=ByxBFsRqYm}.

\bibitem[{Krumm(2008)}]{krumm2008markov}
Krumm J, 2008 \emph{A {M}arkov model for driver turn prediction}. \emph{Society
  of Automotive Engineers (SAE) World Congress, April 2008}.

\bibitem[{Lawler(1985)}]{lawler1985traveling}
Lawler EL, 1985 \emph{The Traveling Salesman Problem: A Guided Tour of
  Combinatorial Optimization} (Wiley-Interscience Series in Discrete
  Mathematics).

\bibitem[{Levenshtein et~al.(1966)}]{levenshtein1966binary}
Levenshtein VI, et~al., 1966 \emph{Binary codes capable of correcting
  deletions, insertions, and reversals} 10:707--710.

\bibitem[{Liao \protect\BIBand{} Liu(2018)}]{8404041}
Liao E, Liu C, 2018 \emph{A hierarchical algorithm based on density peaks
  clustering and ant colony optimization for traveling salesman problem}.
  \emph{IEEE Access} 6:38921--38933.

\bibitem[{Lin et~al.(2014)Lin, Choy, Ho, Chung, \protect\BIBand{}
  Lam}]{Canhong2014}
Lin C, Choy K, Ho G, Chung SH, Lam H, 2014 \emph{Survey of green vehicle
  routing problem: Past and future trends}. \emph{Expert Systems with
  Applications} 41:1118--1138.

\bibitem[{MacGregor \protect\BIBand{} Ormerod(1996)}]{macgregor1996human}
MacGregor JN, Ormerod T, 1996 \emph{Human performance on the traveling salesman
  problem}. \emph{Perception \& Psychophysics} 58(4):527--539.

\bibitem[{Miller, Tucker, \protect\BIBand{} Zemlin(1960)}]{Miller1960IntegerPF}
Miller CE, Tucker AW, Zemlin RA, 1960 \emph{Integer programming formulation of
  traveling salesman problems}. \emph{Journal of the ACM} 7(4):326--329.

\bibitem[{Montoya-Torres et~al.(2015)Montoya-Torres, {López Franco}, {Nieto
  Isaza}, {Felizzola Jiménez}, \protect\BIBand{}
  Herazo-Padilla}]{MONTOYATORRES2015115}
Montoya-Torres JR, {López Franco} J, {Nieto Isaza} S, {Felizzola Jiménez} H,
  Herazo-Padilla N, 2015 \emph{A literature review on the vehicle routing
  problem with multiple depots}. \emph{Computers \& Industrial Engineering}
  79:115--129.

\bibitem[{{OpenStreetMap}(2021)}]{OpenStreetMap}
{OpenStreetMap}, 2021 \emph{{Planet dump retrieved from
  https://planet.osm.org}}. \urlprefix\url{https://www.openstreetmap.org},
  (accessed on 2021-20-11).

\bibitem[{Orman \protect\BIBand{} Williams(2007)}]{Orman2007ASO}
Orman A, Williams HP, 2007 \emph{A survey of different integer programming
  formulations of the travelling salesman problem}. \emph{Advances in
  Computational Management Science}, volume~9, 91--104 (Springer Berlin
  Heidelberg).

\bibitem[{Papadimitriou(1977)}]{papadimitriou1977euclidean}
Papadimitriou CH, 1977 \emph{The {E}uclidean traveling salesman problem is
  {NP}-complete}. \emph{Theoretical Computer Science} 4(3):237--244.

\bibitem[{Pizlo et~al.(2006)Pizlo, Stefanov, Saalweachter, Li, Haxhimusa,
  \protect\BIBand{} Kropatsch}]{pizlo2006traveling}
Pizlo Z, Stefanov E, Saalweachter J, Li Z, Haxhimusa Y, Kropatsch WG, 2006
  \emph{Traveling salesman problem: A foveating pyramid model}. \emph{The
  Journal of Problem Solving} 1(1):8.

\bibitem[{Purkayastha et~al.(2020)Purkayastha, Chakraborty, Saha,
  \protect\BIBand{} Mukhopadhyay}]{Roneeta2020}
Purkayastha R, Chakraborty T, Saha A, Mukhopadhyay D, 2020 \emph{Study and
  analysis of various heuristic algorithms for solving travelling salesman
  problem--a survey}. Mandal JK, Mukhopadhyay S, eds., \emph{Proceedings of the
  Global AI Congress 2019}, 61--70 (Springer Singapore).

\bibitem[{Sengupta, Mariescu-Istodor, \protect\BIBand{}
  Fr{\"a}nti(2018)}]{sengupta2018planning}
Sengupta L, Mariescu-Istodor R, Fr{\"a}nti P, 2018 \emph{Planning your route:
  Where to start?} \emph{Computational Brain \& Behavior} 1(3-4):252--265.

\bibitem[{Speranza \protect\BIBand{} Archetti(2014)}]{Speranza2014}
Speranza M, Archetti C, 2014 \emph{A survey on matheuristics for routing
  problems}. \emph{EURO Journal on Computational Optimization} 2:223--246.

\bibitem[{Taniguchi \protect\BIBand{} Thompson(2002)}]{taniguchi2002modeling}
Taniguchi E, Thompson RG, 2002 \emph{Modeling city logistics}.
  \emph{Transportation Research Record} 1790(1):45--51.

\bibitem[{\text{United Nations}(2018)}]{UN2018}
\text{United Nations}, 2018 \emph{68\% of the world population projected to
  live in urban areas by 2050, says {UN}.}
  \urlprefix\url{https://www.un.org/development/desa/en/news/population/2018-revision-of-world-urbanization-prospects.html},
  (accessed on 2021-20-11).

\bibitem[{Van~Rooij, Stege, \protect\BIBand{} Schactman(2003)}]{van2003convex}
Van~Rooij I, Stege U, Schactman A, 2003 \emph{Convex hull and tour crossings in
  the {E}uclidean traveling salesperson problem: Implications for human
  performance studies}. \emph{Memory \& Cognition} 31(2):215--220.

\bibitem[{Vickers et~al.(2003)Vickers, Lee, Dry, \protect\BIBand{}
  Hughes}]{vickers2003roles}
Vickers D, Lee MD, Dry M, Hughes P, 2003 \emph{The roles of the convex hull and
  the number of potential intersections in performance on visually presented
  traveling salesperson problems}. \emph{Memory \& Cognition} 31(7):1094--1104.

\bibitem[{Wang et~al.(2015)Wang, Ma, Di, Murphey, Qiu, Kristinsson, Meyer,
  Tseng, \protect\BIBand{} Feldkamp}]{wang2015building}
Wang X, Ma Y, Di J, Murphey YL, Qiu S, Kristinsson J, Meyer J, Tseng F,
  Feldkamp T, 2015 \emph{Building efficient probability transition matrix using
  machine learning from big data for personalized route prediction}.
  \emph{Procedia Computer Science} 53:284--291.

\bibitem[{Wiener, Ehbauer, \protect\BIBand{} Mallot(2009)}]{wiener2009planning}
Wiener J, Ehbauer N, Mallot H, 2009 \emph{Planning paths to multiple targets:
  Memory involvement and planning heuristics in spatial problem solving}.
  \emph{Psychological Research PRPF} 73(5):644--658.

\bibitem[{Wiener \protect\BIBand{} Mallot(2003)}]{wiener2003fine}
Wiener JM, Mallot HA, 2003 \emph{'{F}ine-to-coarse' route planning and
  navigation in regionalized environments}. \emph{Spatial cognition and
  computation} 3(4):331--358.

\bibitem[{Ye et~al.(2015)Ye, Wang, Malekian, Lin, \protect\BIBand{}
  Wang}]{ye2015method}
Ye N, Wang Zq, Malekian R, Lin Q, Wang Rc, 2015 \emph{A method for driving
  route predictions based on hidden {M}arkov model}. \emph{Mathematical
  Problems in Engineering} 2015.

\end{thebibliography}


\newpage

\begin{APPENDIX}{Pseudocode}

In this Appendix, we provide the pseudocode to the methodology presented in Section~\ref{sec:model}. Although the pseudocode has been designed on the basis of the Amazon Last-Mile Routing Research Challenge dataset, it can be easily adapted to different data structures. As a side note, the zone ids are not consistent for different stations, therefore, we run our code for each $s \in \mathcal{S}$, with $\mathcal{S}$ the set of stations. The pseudocode used for the learning phase is reported in Algorithm \ref{alg:model_build}, where the function \texttt{ToZoneSeq} converts the realized sequence of stops into a sequence of zones (see Section \ref{build}), while function \texttt{ComputeCountMatrix} computes an asymmetric matrix where the $({i,j})$-th entry represents the number of times a driver went from the $i$-th to the $j$-th zone. Algorithm \ref{alg:model_apply} illustrates the pseudocode used in the prediction phase. Here, in order to compute the cost matrix, as thoroughly described in Section \ref{sec:pred_zsq}, we compute the distance matrix where each entry refers to the travel time between two different zone centers (computed with \texttt{EstimateZoneCenter}).\\
The code used to run the experiments presented in Section~\ref{sec:perf} is open source and will be made available. 

 \begin{algorithm}
\caption{Learning Zone Preferences and Preprocessing} \label{alg:model_build}
\begin{algorithmic}[1]
\Require{$\mathcal{D}_{s}$: route data set corresponding to station $s$}
\Require{\texttt{StopSeqs}: actual sequence of stops for each route in $\mathcal{D}_{s}$}
\Ensure{$\texttt{CountMatrix} \in \mathbb{R}^{(M+1)\times (M+1)}$: counts of station/zone transitions}
\State \texttt{ZoneSequences} $\gets \emptyset$
\State \texttt{AllZones} $\gets \{s\}$ 
\For{\texttt{route} in $\mathcal{D}_{s}$}
    \State \texttt{ZoneSequences[route]} $\gets$ \Call{ToZoneSeq}{\texttt{StopSeqs[route]}}
    \State \texttt{AllZones} $\gets \texttt{AllZones}\cup \{\texttt{ZoneSequences[route]}\}$
\EndFor
\State \texttt{CountMatrix} $\gets$ \Call{ComputeCountMatrix}{\texttt{AllZones}, \texttt{ZoneSequences}}
\end{algorithmic}
\end{algorithm}

\begin{algorithm}
\caption{Predicting Routes}\label{alg:model_apply}
\begin{algorithmic}[1]
\Require{\texttt{CountMatrix}}: output of Algorithm~\ref{alg:model_build}
\Require{\texttt{Zones}, \texttt{Stops}}: zones and stops characterizing the route
\Require{\texttt{station}, \texttt{TravelTimes}}
\Ensure{\texttt{PredictedStops}: predicted sequence of stops}
\State \texttt{ZoneCenters} $\gets$ \Call{EstimateZoneCenter}{\texttt{Zones}, \texttt{Stops}}
\State \texttt{DistanceMatrix} $\gets$ \Call{ComputeDistMatrix}{\texttt{Zones}, \texttt{ZoneCenters}}
\State \texttt{CostMatrix} $\gets$ \Call{ComputeCostMatrix}{\texttt{DistanceMatrix}, \texttt{CountMatrix}}
\State \texttt{PredictedZoneSeq} $\gets$ \Call{TSP}{\texttt{station}, \texttt{Zones}, \texttt{CostMatrix}}
\State \texttt{PredictedStopSeq} $\gets$ [\texttt{station}]
\For{$\texttt{zone}$ in \texttt{PredictedZoneSeq}}
    \State \texttt{PrevStop} $\gets$  \texttt{PredictedStopSeq}[-1] \Comment{Last element of \texttt{PredictedStopSeq}}
    \State \texttt{AddStop} $\gets$ \Call{ClosestToNextZoneCenter}{\texttt{zone}, \texttt{Stops}}
    \State\texttt{PredictedStopSeq} $\gets$ \texttt{PredictedStopSeq} + [\Call{OTSP}{\texttt{PrevStop}, \texttt{Stops}[\texttt{zone}], \texttt{AddStop}, \texttt{TravelTimes}}]
\EndFor
\end{algorithmic}
\end{algorithm}
\end{APPENDIX}

\end{document}